\input  amstex
\input amsppt.sty
\magnification1200
\vsize=23.5truecm
\hsize=16.5truecm
\NoBlackBoxes

\document

\def\crp{\overline{\Bbb R}_+}

\def\rn{{\Bbb R}^n}
\def\rnp{{\Bbb R}^n_+}

\def\rnpm{\Bbb R^n_\pm}
\def\crnp{\overline{\Bbb R}^n_+}

\def\crnpm{\overline{\Bbb R}^n_\pm}
\def\comega{\overline\Omega }
\def\ang#1{\langle {#1} \rangle}
\def\rp{ \Bbb R_+}

\def\Op{\operatorname{Op}}
\def\N{\Bbb N}
\def\R{\Bbb R}
\def\C{\Bbb C}
\def\Z{\Bbb Z}
\def\ol{\overline}
\def\E{\Cal E}
\def\F{\Cal F}
\def\simto{\overset\sim\to\rightarrow}
\def\supp{\operatorname{supp}}

\def\d{d\!@!@!@!@!@!{}^{@!@!\text{\rm--}}\!}

\document
\topmatter
\title
Resolvents for  fractional-order operators with nonhomogeneous local
boundary conditions
\endtitle
\author Gerd Grubb \endauthor
\affil
{Department of Mathematical Sciences, Copenhagen University,\\
Universitetsparken 5, DK-2100 Copenhagen, Denmark.\\
E-mail {\tt grubb\@math.ku.dk}}\endaffil
\rightheadtext{Resolvents}
\abstract
For $2a$-order strongly elliptic operators $P$ generalizing $(-\Delta
)^a$, $0<a<1$, the homogeneous Dirichlet problem
on a bounded open set $\Omega \subset \rn$ has been widely
studied. Pseudodifferential methods have been applied by the present
author when $\Omega $ is smooth; this is extended in a recent joint
work with Helmut Abels showing exact regularity theorems in the scale of $L_q$-Sobolev spaces $H_q^s$
for $1<q<\infty $, when $\Omega $ is $C^{\tau
+1}$ with a finite $\tau >2a$. We now develop this into
existence-and-uniqueness theorems (or Fredholm theorems), by a study of
the $L_p$-Dirichlet realizations of $P$ and
$P^*$, showing that there are finite-dimensional kernels and cokernels lying in $d^aC^\alpha
(\comega)$ with suitable $\alpha >0$, 
$d(x)=\operatorname{dist}(x,\partial\Omega )$. Similar
results are established for $P-\lambda I$, $\lambda \in\C$. The
solution spaces equal $a$-transmission spaces
$H_q^{a(t)}(\comega)$.

Moreover, the results are extended to nonhomogeneous Dirichlet
problems prescribing the local Dirichlet trace
$(u/d^{a-1})|_{\partial\Omega }$. They are solvable in the larger
spaces
 $H_q^{(a-1)(t)}(\comega)$. Furthermore, the nonhomogeneous problem with a
spectral parameter $\lambda \in\C$,
$$ 
    Pu-\lambda u = f \text { in }\Omega ,\quad
     u=0 \text { in }\rn\setminus \Omega ,\quad
     (u/d^{a-1 })|_{\partial\Omega }=\varphi \text{ on }\partial\Omega ,
    $$
is for $q<(1-a)^{-1}$ shown to be uniquely resp.\ Fredholm solvable
when  $\lambda $ is in
the resolvent set resp.\ the spectrum of the $L_2$-Dirichlet
realization.

The results open up for applications of functional analysis
methods. Here we establish solvability results for 
evolution problems with a time-parameter $t$, both in the case of the homogeneous Dirichlet
condition, and the case where a nonhomogeneous Dirichlet trace
$(u(x,t)/d^{a-1}(x))|_{x\in\partial\Omega }$ is prescribed.
\endabstract

\endtopmatter

\subhead 0. Introduction \endsubhead

Fractional-order operators $P$ have been studied extensively in recent
years, the most prominent example being the fractional Laplacian
$(-\Delta )^a$ ($0<a<1$) of order $2a$. They are of interest in
Probability and Finance, as well as in Differential Geometry and
Mathematical Physics.

Let $P$ be a classical pseudodifferential operator of order $2a$,
strongly elliptic with even symbol. From
its definition on $\rn$, one can define its action on open subsets
$\Omega \subset \rn$ in several ways (this is not obvious since $P$ is generally
nonlocal); the most common choice is to let $P$ act on
functions $u$ defined on $\rn$ but vanishing on $\rn\setminus\Omega $
(i.e., supported in $\comega$), and restrict $Pu$ to $\Omega $
afterwards. This leads to the {\it homogeneous restricted fractional Dirichlet problem}$$
Pu=f\text{ on }\Omega ,\quad u=0\text{ on }\rn\setminus\Omega .\tag0.1
$$
 The strategies to study this include methods from potential theory
 and singular integral operator theory, probabilistic methods, and
 pseudodifferential methods.

The present author has worked with pseudodifferential methods \cite{G15a--G19}, leading
to satisfactory results in cases of $C^\infty $-domains $\Omega $ and
operators depending smoothly on $x$. Other methods have allowed far
less smoothness of $\Omega $ (and in some cases of $P$). To breach this
gap, we have in \cite{AG21} with Helmut Abels worked out a theory that
systematically allows  $\Omega $ to be $C^{1+\tau }$ and $P$ to have a $C^\tau $-smooth dependence on $x$,
for finite positive $\tau $, leading to results in Sobolev-type
spaces $H_q^s$ (Bessel-potential spaces) with parameter $s$ limited by $\tau $, and with
corollaries in H\"older spaces.  We work in this paper with such
operators and domains, under the basic hypothesis that $\tau
>2a$; this is replaced by $\tau >2a+1$ if a nonhomogeneous boundary condition enters. 

An important point in the present investigations is to enhance the regularity
results of \cite{AG21} with genuine solvability results: Theorems about
existence and uniqueness of solutions, or, when relevant, Fredholm
solvability. 
The point of departure is here the properties of the $L_2$-Dirichlet
realization $P_{D,2}$ of $P$ defined variationally from the
sesquilinear form $\int_{\Omega }Pu\, \bar v\, dx$ on $\dot
H^a(\comega)$ (the functions in the Sobolev space $H^a(\rn)$ supported
in $\comega$). By compact embeddings, $P_{D,2}$ 
has a discrete spectrum $\Sigma \subset \C$ consisting of
eigenvalues $\lambda $ with finite dimensional eigenspaces $N_\lambda
$. 
The following issues will be addressed:
\smallskip

\noindent {\bf 1)} The {\it regularity of eigenfunctions} $u_\lambda $, i.e., nontrivial solutions of 
$$
Pu_\lambda =\lambda u_\lambda \text{ on }\Omega ,\quad u_\lambda =0\text{ on }\rn\setminus\Omega ,\tag0.2
$$
 $\lambda \in\C$.
It is shown that the possible $\lambda $-values belong to  $\Sigma $,
and the eigenfunctions 
lie in
$d^aC^\alpha  (\comega)$ for suitable $\alpha >0$;  $d=\operatorname{dist}(x,\partial\Omega )$. See Theorem 4.4 and Corollary 4.5 below. The
regularity of eigenfunctions has been studied earlier for $(-\Delta
)^a$ e.g.\ in Chen and Song \cite{CS05},  Servadei and Valdinoci
\cite{SV13} and Ros-Oton and Serra \cite{RS14},  for smooth $P$ and
$\Omega $  in Grubb \cite{G15b}, and for $\Omega $
 equal to a ball in  Dyda,
 Kuznetsov and Kwasnicki \cite{DKK17}. We also describe the
 eigenfunctions for $P^*$, in Theorem 4.15.
\smallskip

\noindent {\bf 2)} The {\it structure of the $L_q$-Dirichlet realization}
$P_{D,q}$   of $P$ for 
functions satisfying (0.1), $1<q<\infty $. In particular an investigation of $P^*$
and a proof that  $P_{D,q}$ in $L_q(\Omega )$ and $(P^*)_{D,q'}$ in
$L_{q'}(\Omega )$ are adjoints ($\frac1q + \frac1{q'}=1$), and are
Fredholm operators with the same kernel and cokernel as in the case
$q=2$. See Theorem 4.16 $1^\circ$ below. We are not aware of other
studies of $P_{D,q}$ for $q\ne 2$ in nonsmooth cases
with the precision that $D(P_{D,q})$ equals the $a$-transmission space $H_q^{a(2a)}(\comega)$.
\smallskip

\noindent {\bf 3)} A treatment of the {\it resolvent
problem} for $P_{D,q}$:
 $$
    (P-\lambda )u = f \text { in }\Omega ,\quad     u =0 \text { in }\rn\setminus \Omega , \tag0.3
 $$
showing unique solvability when $\lambda \in \C\setminus\Sigma $  and Fredholm
solvability when $\lambda \in\Sigma $, in appropriate Sobolev-type function
spaces. See Theorem 4.16 $2^\circ$ and Theorem 4.17
 below. Also results with a higher regularity parameter $s>0$,
 including H\"older spaces, are
 obtained, see Corollary 4.9 and Theorem 4.18. 
 
\smallskip
\noindent {\bf 4)} A treatment of the {\it local nonhomogeneous Dirichlet
problem} for $P$:
 $$
 \aligned
    Pu& = f \text { in }\Omega ,\\
     u&=0 \text { in }\rn\setminus \Omega ,\\
     (u/d^{a-1 })|_{\partial\Omega }&=\varphi \text{ on }\partial\Omega ,
     \endaligned \tag0.4
 $$
showing regularity results (Theorem 3.7 and Corollary 3.9), and unique or Fredholm solvability (Theorem 5.1 and
Corollary 5.2), in  $H_q^s$-related function
spaces  and H\"older spaces.
\smallskip

\noindent {\bf 5)} A treatment  of {\it local nonhomogeneous Dirichlet problems with a spectral
parameter} $\lambda \in\C $:
 $$
 \aligned
    Pu-\lambda u& = f \text { in }\Omega ,\\
     u&=0 \text { in }\rn\setminus \Omega ,\\
     (u/d^{a-1 })|_{\partial\Omega }&=\varphi \text{ on }\partial\Omega ,
     \endaligned \tag0.5
 $$
when $q<(1-a)^{-1}$. There is  unique solvability
when $\lambda \in\C\setminus\Sigma$, Fredholm solvability  when $\lambda \in\Sigma $. See Theorem 5.4
below. The study of such problems  was initiated by Chan, Gomez-Castro and
Vazquez \cite{CGV21} (with a more implicit
formulation of
the boundary condition, and $P=(-\Delta )^a$ acting in weighted $L_1$-spaces). 
\smallskip

\noindent {\bf 6)} Solvability of {\it evolution problems}
$$\aligned
Pu(x,t)+\partial_tu(x,t)&=f(x,t)\text{ on }\Omega \times I ,
\quad I=\,]0,T[\,,\\
u(x,t)&=0\text{ on }(\Bbb R^n\setminus\Omega )\times I, \\ 
u(x,0)&=0;
\endaligned \tag0.6
$$
possibly with a nonhomogeneous boundary condition
$$
(u(x,t)/d^{a-1}(x))|_{x\in\partial\Omega }=\psi(x,t) \text{ on }\partial\Omega \times I.\tag0.7
$$
See Theorems
6.2--6.5 below. The results without condition (0.7) are a
straightforward extension of results for smooth cases shown in
\cite{G18a}, \cite{G18b}; there are earlier results for
$x$-independent operators on nonsmooth domains with H\"older
estimates by  Ros-Oton with Fernandez-Real and Vivas
\cite{FR17}, \cite{RV18}. Evolution problems prescribing
the nonhomogeneous boundary condition (0.7) have to our knowledge not been studied before. 
\smallskip

A large part of the results are new even for $(-\Delta )^a$.

It is a pervading fact in all these results that the exact operator domains are found; they
have the form of $a$-transmission spaces $H_q^{a(s+2a)}(\comega)$ in
cases with homogenous Dirichlet condition, and $(a-1)$-transmission
spaces $H_q^{(a-1)(s+2a)}(\comega)$ in cases with  nonhomogeneous
Dirichlet condition.

It should be noted that there exist several interpretations of what a
nonhomogeneous Dirichlet condition could be. A frequently studied
possibility is to prescribe an
exterior value of $u$,
$$
u=g \text{ on }\rn \setminus\Omega ;\tag0.8
$$
then the problem can be reduced to the homogeneous Dirichlet problem by
subtraction of a suitable extension of $g$ to $\Omega $ (as described e.g.\ in
\cite{G14b}). Problems with the condition (0.8) are global, involving
all of $\rn$. Our choice of nonhomogeneous Dirichlet condition is
$$
(u/d^{a-1})|_{\partial\Omega }=\varphi \text{ on }\partial\Omega ;\tag0.9
$$
it is localized to $\partial\Omega $, even pointwise. In Section 1 below,
we explain by comparison with $\Delta $ why this choice is natural for $(-\Delta )^a$. 

\bigskip

\noindent {\it Plan of the paper:} Section 1 introduces the local
nonhomogeneous Dirichlet condition for $2a$-order operators. Section 2
sets up the terminology, introducing function spaces,
pseudodifferential operators, and the special $\mu $-transmission
spaces and their role in the definition of weighted boundary
values. Section 3 recalls the regularity result for the homogeneous
Dirichlet problem known from \cite{AG21}, and establishes regularity
results for the nonhomogeneous Dirichlet problem. In Section 4,
realizations $P_{D,q}$ of the homogeneous Dirichlet problem for $P$ in
$L_q(\Omega )$ ($1<q<\infty $) are studied. Regularity of
eigenfunctions is shown, also for $P^*$, and the resolvent
$(P_{D,q}-\lambda )^{-1}$ is set up for $\lambda \notin\Sigma $, where
$\Sigma $ denotes the spectrum of $P_{D,2}$. Fredholm properties
are established for $P_{D,q}-\lambda $ when $\lambda \in\Sigma
$. Section 5 shows the unique or Fredholm solvability of the
nonhomogeneous Dirichlet problem, and shows how $P$ can be replaced by
$P-\lambda $ when $q<(a-1)^{-1}$. Finally, Section 6 treats evolution
problems in cases where there is a uniform  norm estimate of the resolvent.

\subhead 1. A simple introduction to the local
nonhomogeneous Dirichlet problem for $(-\Delta )^a$ \endsubhead

\subsubhead 1.1 Standard elliptic boundary problems \endsubsubhead
First recall some facts about the Laplacian $A=-\Delta $ (they are
also true for strongly elliptic
second-order differential operators $A$ with smooth coefficients,
satisfying $\operatorname{Re}(Au,u)>0$ when $u\in C_0^\infty (\Omega
)\setminus \{0\}$).

Let $\Omega \subset \rn$ be a bounded smooth domain. For simplicity,
consider solvability just in a $C^\infty $ setting. For any
 $\mu >-1$ we define the space $\E_\mu $ by: 
$$
\E_\mu (\comega)=e^+d_0^\mu C^\infty (\comega),\tag1.1
$$
where $d_0(x)$ is a function equal to $\operatorname{dist}(x,\partial\Omega )$
on a neighborhood of $\Omega $, extended as a positive $C^\infty $-function
to the rest of $\Omega $. ($d_0$ can be replaced by an equivalent
function $d$, see (2.2)ff.\ below.)
Here  $e^+$ denotes extension by zero on
$\rn\setminus\Omega $; it is relevant in the consideration of nonlocal
operators. We shall also use the notation $r^+$ that indicates
restriction from $\rn$ to $\Omega $.

The functions in $\E_0 $ have
Taylor expansions in the neighborhood of each boundary point, where
$x$ equals $(x',x_n)$ in local coordinates $(x',x_n)\in \R^{n-1}\times\rp$:
$$
u(x)=u_0(x')+u_1(x')x_n+\tfrac12 u_2(x')x_n^2+\dots\text{ for }x_n>0;\tag1.2
$$
here $u_0=\gamma _0u=\lim_{x_n\to 0+}u(x',x_n)$, and $u_j=\gamma _ju=\gamma
_0(\partial_{x_n}^ju)$, for $j= 1,2,\dots$, with 
the usual notation for boundary values.

For integer values $k\ge 1$, the functions in $\E_k(\comega)$ have Taylor
expansions at the boundary, like (1.2) but skipping the first $k$
terms. For example for $u\in \E_1(\comega)$,
$$
u(x)=u_1(x')x_n+\tfrac12 u_2(x')x_n^2+\dots\text{ for }x_n>0.\tag1.3
$$

It is well-known that the {\it nonhomogeneous Dirichlet problem} for $A$:
$$
Au=f \text{ on }\Omega , \quad \gamma _0u=\varphi \text{ on }\partial\Omega ,\tag1.4
$$
is uniquely solvable for any $f\in C^\infty (\comega)$,  $\varphi \in
C^\infty (\partial\Omega )$, {\bf with} $u\in \E_0(\comega)\simeq C^\infty (\comega)$. 

In particular, the {\it homogeneous Dirichlet problem} for $A$:
$$
Au=f \text{ on }\Omega , \quad \gamma _0u=0 \text{ on }\partial\Omega ,\tag1.5
$$
is uniquely solvable for any $f\in C^\infty (\comega)$; {\bf here}
$u\in \E_1(\comega)$, {\bf satisfying }(1.3). 

\subsubhead 1.2 The fractional Laplacian \endsubsubhead
Now consider $P=(-\Delta )^a$ with $0<a<1$. It is of order $2a$, and the
symbol $p(\xi )=|\xi |^{2a}$ is {\it even} in $\xi $ (satisfies $p(-\xi
)=p(\xi )$). It was proved in \cite{H66,H85} (and presented in detail
in \cite{G15a}) by a fine analysis of what
happens at the boundary, that $r^+P$ has a good meaning on
$\E_{a+k}(\comega)$ with $k$ integer $\ge -1$, mapping
$$
r^+P\colon \E_{a+k}(\comega)\to C^\infty (\comega).\tag1.6
$$

The {\it homogeneous restricted Dirichlet problem} for $P$ is generally
agreed to be the problem
$$
Pu=f \text{ on }\Omega , \quad  u=0 \text{ on }\rn\setminus\comega .\tag1.7
$$
It is well-known (by a variational argument) 
that this problem is uniquely solvable
when $u$ is a priori sought in $\dot H^a(\comega)=\{u\in H^a(\rn)\mid
u=0 \text{ on }\rn\setminus\Omega \}$ and $f$ is given in
$L_2(\Omega )$. The regularity question is about how a higher
regularity of $f$ implies a higher regularity of $u$.

It is shown in  \cite{G15a} that when $f\in
C^\infty (\comega)$, the solution of (1.7) is
in fact in $\E_a(\comega)$. Thus $r^+P$ defines a homeomorphism:
$$
r^+P\colon \E_a(\comega)\simto C^\infty (\comega).\tag1.8
$$

Note that, by multiplication by $d^a$, one has from (1.2) in local coordinates:
$$
\text{when }u\in \E_a,\quad
u(x)=v_0(x')x_n^a+v_1(x')x_n^{a+1}+\tfrac12
v_2(x')x_n^{a+2}+\dots\text{ for }x_n>0,\tag1.9
$$
where $v_0=\gamma _0(u/x_n^a)$, $v_1=\gamma _1(u/x_n^a)$, etc. Then
$\E_a$ has a role parallel to that of $\E_1$ in the standard
homogeneous Dirichlet problem (1.5), cf.\ (1.3).

Analogously to the nonhomogeneous standard Dirichlet problem (1.4) we now
consider $\E_{a-1}$, which will have a role parallel to that of $\E_0$ in
the following {\bf nonhomogeneous local Dirichlet problem for }$P$:
$$
Pu=f \text{ on }\Omega , \quad  u=0 \text{ on }\rn\setminus\comega,\quad
\gamma _0(u/d^{a-1})=\varphi  ,\tag1.10
$$
where $u$ is sought in $\E_{a-1}$.
Indeed, the boundary behavior of functions in $\E_{a-1}$ is 
$$
\text{when }u\in \E_{a-1},\quad
u(x)=w_0(x')x_n^{a-1}+w_1(x')x_n^{a}+\tfrac12
w_2(x')x_n^{a+1}+\dots\text{ for }x_n>0,\tag1.11
$$
where  $w_0=\gamma _0(u/x_n^{a-1})$, $w_1=\gamma _1(u/x_n^{a-1})$, etc.
Note that the expansion is similar to that in (1.9), the only difference
being that the coefficient $w_0$
vanishes there, i.e.,
$$
\E_a\text{ is the subset of }\E_{a-1} \text{ where }\gamma
_0(u/x_n^{a-1})=0.\tag1.12 
$$
Using that a given $\varphi  \in C^\infty (\partial\Omega )$ can be lifted
to a function $z\in \E_{a-1}$ such that  $\gamma
_0(z/d^{a-1})=\varphi  $ (namely, locally, $z(x)=\varphi (x')x_n^{a-1}$), we
get immediately the unique solvability of the nonhomogeneous problem
(1.10) from the solvability of the homogeneous problem  (1.7).

For generalizations $P$ of the fractional Laplacian with smooth, even
symbol, one finds the same results. This shows the interesting fact that $\E_a$ is
{\it universal} as the solution space for the homogeneous Dirichlet
problem, and that $\E_{a-1}$ similarly plays a universal role for our
nonhomogeneous  Dirichlet problem.

Let us mention briefly that one can also define a {\it local Neumann condition} for $P$ in analogy
with the standard Neumann condition for $A$, namely by prescribing
$\gamma _1(u/d^{a-1})$. Also here there are general solvability
results; more details are found in e.g.\ \cite{G14b}, \cite{G18}.

Note that the functions  $u$ in $ \E_{a-1}$ {\bf blow up} like
$d^{a-1}$ at the
boundary at the points where
$\gamma _0(u/d^{a-1})$ does not vanish; this is a natural fact in the theory.

\subsubhead 1.3 Results in Sobolev spaces \endsubsubhead
There is now the question of how these problems are treated in more
general function spaces. For example, in terms of $L_2$-Sobolev
spaces, the homogeneous  Dirichlet problem for the
Laplacian (1.5)  is solved in $H^2(\Omega )\cap H^1_0(\Omega )$ when
$f\in L_2(\Omega )$,
and the nonhomogeneous  Dirichlet problem for the
Laplacian (1.4)  is solved in $H^2(\Omega )$ when
$f\in L_2(\Omega )$, $\varphi \in H^{\frac32}(\partial\Omega )$. The
same spaces enter when $\Delta $ is replaced by a strongly elliptic
second-order differential operator $A$ with smooth coefficients.

There are corresponding results for $(-\Delta )^a$:
The homogeneous  Dirichlet problem (1.7)  is solved in $H^{a(2a)}(\comega )$ when
$f\in L_2(\Omega )$,
and the nonhomogeneous  Dirichlet problem (1.10)  is solved in $H^{(a-1)(2a)}(\comega )$ when
$f\in L_2(\Omega )$, $\varphi \in H^{a+\frac12}(\partial\Omega )$.
Here the so-called transmission spaces $H^{a(s)}(\comega)$ and
$H^{(a-1)(s)}(\comega)$ enter; they were defined in \cite{G15a}
(building on \cite{H66}), and are important since they give exact
information. $H^{a(2a)}(\comega)$ takes the place of $H^2(\Omega )\cap
H^1_0(\Omega )$ and contains the space $\E_a(\comega)$, whereas
$H^{(a-1)(2a)}(\comega)$ takes the place of $H^2(\Omega )$ and contains $\E_{a-1}(\comega)$.
(Their definition is recalled in the general preliminaries section below.)
Also here the spaces are universal; the same spaces enter when $(-\Delta )^a$ is replaced by a
strongly elliptic 
pseudodifferential operator $P$ of order $2a$ with even symbol.

\example{Remark 1.1}
The nonhomogeneous Dirichlet problem (1.10) for $(-\Delta )^{a}$ was
proposed simultaneously and independently in the works \cite{G15a} and
Abatangelo \cite{A15}, in two very different formulations. Ours was in
the style explained above, whereas Abatangelo formulated the problem
in its relation to a Green's function and a 
representation of the solution by a sum of integrals over $\Omega $
and $\partial\Omega $. In \cite{A15}, the boundary value is
somewhat implicitly formulated as a term $Eu$ defined by integrals, and only appears in the form $c\gamma
_0(u/d^{a-1})$ in the case where $\Omega $ is a ball. The
word ``large solution'' is introduced to underline the blow-up (like
$d^{a-1}$) that solutions with nonzero continuous boundary data will have at the boundary.

It has been known to many people as an accepted fact (or folklore)
that prescribing $Eu$ is
equivalent  to prescribing $\gamma _0(u/d^{a-1})$. A proof
that $Eu$
is proportional to  $\gamma _0(u/d^{a-1})$  for the
fractional Laplacian is included as  Appendix
A.1 below.
A different proof is given in
App.\ B of \cite{CGV21}, with another
proportionality factor (see Theorem A.4ff.\ below).
\endexample

\subhead 2. Preliminaries \endsubhead

\subsubhead 2.1 Function spaces \endsubsubhead

The space $C^k(\rn)\equiv C^k_b(\rn)$ consists  of $k$-times differentiable
functions with uniform norms $\|u\|_{C^k}=\sup_{|\alpha |\le k,x\in\rn}|D^\alpha
u(x)|$ ($k\in{\Bbb N}_0$), and the H\"older
spaces $C^\tau   (\rn)$, $\tau  =k+\sigma $ with $k\in{\Bbb N}_0$,
$0<\sigma <1$, also denoted $C^{k,\sigma } (\rn)$, consists of
function  with norms
$\|u\|_{C^\tau  }=\|u\|_{C^k}+\sup_{|\alpha |= k,x\ne y}|D^\alpha
u(x)-D^\alpha u(y)|/|x-y|^\sigma $. The latter definition extends to
Lipschitz spaces $C^{k,1} (\rn)$. There are similar spaces over subsets
of $\rn$. We denote $C^\infty_b(\rn)= \bigcap_{k\in\N} C^k_b(\rn)$.

The halfspaces $\rnpm$ are defined by
 $\rnpm=\{x\in
{\Bbb R}^n\mid x_n\gtrless 0\}$, with points denoted  $x=(x',x_n)$,
$x'=(x_1,\dots, x_{n-1})$. For a given real function $\zeta \in C^{1+\tau}
(\R^{n-1})$ (some $\tau >0$), we define the curved halfspace
$\rn_\zeta  $ by $$
\rn_\zeta = \{x\in\R^n\mid x_n>\zeta (x')\};\tag2.1
$$
it is
 a $C^{1+\tau }$-domain. (The function $\zeta $ was denoted $\gamma $ in
 \cite{AG21}; we change the name to avoid confusion with the notation for trace
 operators $\gamma _j$.)

By a bounded $C^{1+\tau }$-domain $\Omega $  we mean the following:
 $\Omega \subset\rn$ is open and bounded, and every boundary point
$x_0$ has an open neighborhood $U$ such that, after a translation of
$x_0$ to $0$ and a suitable rotation, $U\cap \Omega $ equals $U\cap \rn_\zeta 
$ for a function $\zeta  \in C^{1+\tau }(\R^{n-1})$ with $\gamma
(0)=0$.

Restriction from $\R^n$ to $\rnpm$ (or from
${\Bbb R}^n$ to $\Omega $ resp.\ $\complement\comega= \rn \setminus \comega$) is denoted $r^\pm$,
 extension by zero from $\rnpm$ to $\R^n$ (or from $\Omega $ resp.\
 $\complement\comega$ to ${\Bbb R}^n$) is denoted $e^\pm$. (The
 notation is also used for $\Omega =\rn_\zeta  )$.) Restriction
 from $\crnp$ or $\comega$ to $\partial\rnp$ resp.\ $\partial\Omega $
 is denoted $\gamma _0$.

When $\Omega $ is a $C^{1+\tau }$-domain, we denote by $d(x)$ 
 (as in \cite{G15, Def.\ 2.1} for the $C^\infty
$-case) a function  that is $C^{1+\tau } $ on $\comega$,
positive on $\Omega $ and vanishes only to the first order on
$\partial\Omega $ (i.e., $d(x)=0$ and $\nabla d(x)\neq 0$ for $x\in
\partial\Omega$).  On bounded sets it satisfies near $\partial\Omega $:
$$
  C^{-1}d_0(x)\le d(x)\le Cd_0(x)\tag 2.2
$$
with $C>0$, where $d_0(x)$ equals
$\operatorname{dist}(x,\partial\Omega )$ on a  neighborhood of
$\partial\Omega $ and is extended as a correspondingly smooth positive function on
$\Omega $.
When $\tau \ge 1$,
$d_0$ itself can be taken $C^{1+\tau}$ (as explained e.g.\ in
\cite{AG21}),
then moreover, 
 $d/d_0$ is a positive $C^{\tau }$-function on $\comega$.
   
We take $d_0(x)=x_n$ in the
case of $\rnp$. For $\rn_\zeta $, the function  $d(x)=x_n-\zeta
(x')$ satisfies (2.2) when 
the extension of $d_0(x)$ is suitably chosen for large $x_n$ (cf.\ e.g.\ \cite{AG21}).

The Bessel-potential spaces
$H^s_q({\Bbb R}^n)$ are defined
for $s\in{\Bbb R}$, $1<q<\infty $, by 
$$
H_q^s(\R^n)=\{u\in \Cal S'({\Bbb R}^n)\mid \F^{-1}(\ang{\xi }^s\hat u)\in
L_q(\R^n)\},\tag2.3
$$
where $\Cal F$ is the Fourier transform  $\hat
u(\xi )=\Cal F
u(\xi )= \int_{{\Bbb R}^n}e^{-ix\cdot \xi }u(x)\, dx$, and the
function $\ang\xi $ equals $(|\xi |^2+1)^{\frac12}$.  For  $q=2$, this is the scale of Sobolev spaces,
where the index $2$ is usually omitted. $\Cal S'(\rn)$ is the Schwartz
space of temperate distributions, the dual space of $\Cal S(\rn)$
(the space of rapidly decreasing $C^\infty $-functions).

For $s\in {\Bbb
N}_0=\{0,1,2,\dots\}$, the spaces $H_q^s(\rn)$ are also denoted $W_q^s({\Bbb R}^n)$ or $W^{s,q}({\Bbb R}^n)$
in the literature. We moreover need to refer to the Besov
spaces  $B^s_{q,q}(\rn)$, also denoted $B^s_q(\rn)$,
that coincide with the
$W^s_q$-spaces when $s\in \rp\setminus \N$.
They necessarily enter in connection with
 boundary value problems in an $H^s_q$-context,
 because they are the correct range spaces for
trace maps $\gamma _ju=(\partial_n^ju)|_{x_n=0}$:
$$
\gamma _j\colon \ol H^s_q(\rnp), \ol B^s_q(\rnp) \to
B_q^{s-j-\frac1q}({\Bbb R}^{n-1}), \text{ for }s-j-\tfrac1q >0,\tag2.4
$$
 (cf.\ (2.5)), surjectively and with a continuous right inverse; see e.g.\ the overview in
the introduction to \cite{G90}. For $q=2$, the two scales $H_q^s$ and
$B_q^s$ are
identical, but for $q\ne 2$ they are related by strict inclusions: $
H^s_q\subset B^s_q\text{ when }q>2$, $H^s_q\supset B^s_q\text{ when
}q<2$.

Along with the spaces $H^s_q({\Bbb R}^n)$ defined in (2.3), there are
the two scales of spaces associated with $\Omega $ for $s\in{\Bbb R}$:
$$
\aligned
\ol H_q^{s}(\Omega)&=\{u\in \Cal D'(\Omega )\mid u=r^+U \text{ for some }U\in
H_q^{s}(\R^n)\}, \text{ the {\it restricted} space},\\
\dot H_q^{s}(\comega)&=\{u\in H_q^{s}({\Bbb R}^n)\mid \supp u\subset
\comega \},\text{ the {\it supported} space;}
  \endaligned\tag2.5
$$
here $\operatorname{supp}u$ denotes the support of $u$ (the complement
 of the largest open set where $u=0$). 
 $\ol H_q^s(\Omega )$ is in other texts often denoted  $H_q^s(\Omega )$  or
$H_q^s(\comega )$, and $\dot H_q^{s}(\comega)$ may be indicated with a
ring, zero or twiddle;
the current notation stems from H\"ormander \cite{H85,
 App. B.2}.
There is an identification of $\ol H_q^s(\Omega )$ with the dual space
of $\dot H_{q'}^{-s}(\comega)$, $\frac1{q'}=1-\frac1q\,$, in terms of a
duality extending the sesquilinear scalar product
$\ang{f,g}=\int_\Omega f\,\ol g\, dx$.

Besides for the $H^s_q$ and $B^s_q$-spaces, there are in \cite{G14b} for
$C^\infty $-domains 
established the relevant results in many other scales of spaces,
namely Besov spaces $B^s_{p,q}$ for $1\le p,q \le \infty $
and Triebel-Lizorkin spaces $F^s_{p,q}$ (for the same $p,q$ but with
$p<\infty $). Here we just want to mention
the H\"older-Zygmund scale $B^s_{\infty ,\infty }$, also denoted
$C^s_*$. The space $C^s_*$ identifies with the
H\"older space $C^s$ when $s\in
\rp\setminus {\Bbb N}$, and for positive integer $k$ satisfies
$ C^{k-\varepsilon }\supset C^k_*\supset 
C^{k-1,1}\supset  C^k_b$ for small $\varepsilon >0$; moreover,
$C^0_*\supset L_\infty \supset C^0_b$ (with strict inclusions everywhere). Similarly to (2.5), we denote the
spaces of restricted, resp.\ supported elements
$$
\aligned
\ol C_*^{s}(\Omega)&=\{u\in \Cal D'(\Omega )\mid u=r^+U \text{ for some }U\in
C_*^{s}(\R^n)\},\\
\dot C_*^{s}(\comega)&=\{u\in C_*^{s}({\Bbb R}^n)\mid \supp u\subset
\comega \}.
\endaligned\tag2.6
$$
The star can be omitted when $s\in \rp\setminus \N$ (then we shall
 often write $\ol C^s(\Omega )$ in the more established notation
 $C^s(\comega)$). H\"older
 spaces  over $C^{1+\tau }$-domains $\Omega $ are used in \cite{AG21}.

\subsubhead 2.2 Pseudodifferential operators  \endsubsubhead

A {\it pseudodifferential operator} ($\psi $do) $P$ on ${\Bbb R}^n$ is
defined from a function $p(x,\xi )$  on ${\Bbb
R}^n\times{\Bbb R}^n$, called the {\it symbol},  by 
$$
Pu=\Op (p(x,\xi ))u 
=(2\pi )^{-n}\int_{\rn} e^{ix\cdot\xi
}p(x,\xi )\hat u(\xi)\, \d\xi =\Cal F^{-1}_{\xi \to x}(p(x,\xi )\F u(\xi
)),\tag2.7
$$
using the Fourier transform $\F$. An introduction to $\psi $do's is
given e.g.\ in \cite{G09, Ch.\ 7--8}. A description with more references
and an inclusion of results for operators with nonsmooth symbols can
be found in
\cite{AG21}. We shall here just give a quick summary of definitions and
consequences that we need in the present paper.

The space $S^m_{1,0}({\Bbb R}^n\times{\Bbb R}^n)$ of symbols
$p$ of order $m\in{\Bbb R}$ consists of the complex
$C^\infty $-functions $p(x,\xi )$
such that $\partial_x^\beta \partial_\xi ^\alpha p(x,\xi
)$ is $O(\ang\xi ^{m-|\alpha |})$ for all $\alpha ,\beta $, for some
$m\in{\Bbb R}$, with global estimates in $x\in{\Bbb R}^n$.
$P$ is then of order $m$. It  maps $H^s_q({\Bbb R}^n)$ continuously into
$H^{s-m}_q ({\Bbb R}^n)$ for all $s\in{\Bbb R}$.

$P $ with symbol $p\in  S^m_{1,0}({\Bbb R}^n\times{\Bbb R}^n)$ is said to be {\it classical} when 
$p$  
has an asymptotic expansion $p(x,\xi )\sim \sum_{j\in{\Bbb
N}_0}p_j(x,\xi )$ with $p_j$ homogeneous in $\xi $ of degree $m-j$ for
all $|\xi |\ge 1$ and $j\in\N_0$, such that
$$
\partial_x^\beta \partial_\xi ^\alpha \bigl(p(x,\xi )-
{\sum}_{j<J}p_j(x,\xi )\bigr) \text{ is }O(\ang\xi ^{m-\alpha -J})\text{ for
all }\alpha ,\beta \in{\Bbb N}_0^n, J\in{\Bbb N}_0.  \tag2.8
$$
The space of classical symbols is denoted $S^m({\Bbb R}^n\times{\Bbb R}^n)$.
 For a 
complete theory one adds to these operators 
the {\it smoothing operators} (mapping any  $H^s_q({\Bbb R}^n)$ into
$\bigcap_tH^t_q ({\Bbb R}^n)$), regarded as operators of order
$-\infty $. (For example, $(-\Delta )^a$ fits into the calculus when
it is 
written as $\Op((1-\eta (\xi ))|\xi |^{2a})+\Op(\eta (\xi
)|\xi |^{2a})$, where $\eta (\xi )$ is a $C^\infty $-function that
equals $1$ for $|\xi |\le \frac12$ and 0 for $|\xi |\ge 1$; the 
second term is smoothing.)

Symbols with finite smoothness in $x$ are defined as follows: The
symbol space \linebreak $C^\tau S^m_{1,0}(\R^{n}\times \rn)$ for $\tau>0$, $m\in\R$, consists of functions
$p \colon  \R^{n}\times \R^n\to \C$ that are continuous w.r.t.\
$(x,\xi)\in  \R^{n} \times \R^n$ and $C^\infty $ with respect to $\xi\in
\R^n$, such that for every $\alpha\in\N_0^n$  we have:
$\partial_\xi^\alpha p(x,\xi)$ is in $C^\tau (\R^{n})$ with respect to
$x$ and satisfies for all $\xi\in\R^n$, $\alpha \in {\Bbb N}_0^n$, 
$$
  \|\partial_\xi^\alpha p(\cdot,\xi)\|_{C^\tau  (\R^{n})}\leq C_\alpha
  \ang{\xi}^{m-|\alpha|} ,
\tag2.9
$$
with $C_\alpha>0$. The symbol space is a Fr\'echet space with the semi-norms
$$
  |p|_{k,C^\tau S^m_{1,0}(\R^{n}\times \R^n)}:= \max_{|\alpha|\leq k} \sup_{\xi\in\R^n} \ang{\xi}^{-m+|\alpha|}\|\partial_\xi^\alpha p(\cdot,\xi)\|_{C^\tau  (\R^{{n}})}\quad \text{for }k\in\N_0.\tag2.10
$$
For such symbols there holds when $\tau >0$:
$$
    \Op(p)\colon H^{s+m}_q (\R^n)\to H^s_q(\R^n)\quad \text{for all }|s|<\tau,\tag2.11
$$
where the operator norm for each $s$ 
is estimated by a finite system of symbol seminorms (depending on $s$).

 As
explained in detail in \cite{G14a, Sect.\ 2.3}, the operators can be approximated by operators with
 smooth symbols:
When $p\in  C^\tau
S^{m}_{1,0}(\rn\times\rn)$, it is approximated in the seminorms of 
$  C^{\tau '}
S^{m}_{1,0}(\rn\times\rn)$, any $\tau '<\tau $, by the convolutions in $x$
with an approximate unit: $\varrho _k(x)=k^n\varrho (kx)$ for a
$\varrho \in C_0^\infty (\rn)$ with $\|\varrho \|_{L_1}=1$; here
$p_k=\varrho _k*p\in S^{m}_{1,0}(\rn\times\rn)$.
Hence, taking $\tau '>|s|$, $P_k=\Op(\varrho _k*p)$, 
$$
\|P-P_k\|_{\Cal L(H_q^{s+m}(\rn),H_q^{s}(\rn))}\to 0 \text{ for }k\to \infty
,\text{ when }|s|< \tau '=\tau -\varepsilon  ,\tag2.12
 $$
where $\varepsilon >0$ can be taken arbitrarily small.

The subspace of {\it classical} symbols $C^\tau
S^m(\R^{n}\times \rn)$ consists of those functions that moreover have
expansions into terms  $p_j$ homogeneous in $\xi $ of degree $m-j$ for
$|\xi |\ge 1$, all $j$, such that for all $\xi\in\R^n$, $\alpha
\in{\Bbb N}_0^n$, $ J\in{\Bbb N}_0$,
$$
\|\partial_\xi ^\alpha \bigl(p(\cdot,\xi )-
{\sum}_{j<J}p_j(\cdot,\xi )\bigr)\|_{C^\tau  (\R^{n})}\leq C_{\alpha
,J}\ang{\xi}^{m-J-|\alpha|}.
\tag2.13
$$

A classical symbol $p(x,\xi )$ (and the associated operator $P$) is said to be   {\it strongly elliptic}
when $\operatorname{Re}p_0(x,\xi )\ge c|\xi |^m $ for $|\xi |\ge 1$,
with $c>0$.
Moreover, a classical $\psi $do $P=\Op(p(x,\xi ))$ of
order $m\in\R$ is said
to be 
{\it even}, when the terms in the symbol expansion $p\sim\sum_{j\in\N_0}p_j$ satisfy
$$
p_j(x,-\xi )=(-1)^jp_j(x,\xi )\quad \text{ for all }x\in \rn,|\xi |\ge 1, j\in \Bbb N_0.
\tag2.14
$$
(The word ``even'' is short for {\it even-to-even
  parity}, meaning that  the
terms with even $j$ are even in $\xi $,  the
terms with odd $j$ are odd in $\xi $.)

\subsubhead 2.3 $\mu $-transmission
spaces \endsubsubhead

The following is a rapid introduction to $\mu $-transmission
spaces, which were presented in full detail
in \cite{G15a}, and extended to nonsmooth
domains $\Omega $ in \cite{AG21}.

For some types of $\psi $do's, namely those of integer order belonging
to the Boutet de Monvel calculus (as initiated in \cite{B71}, see
e.g.\ \cite{G90}, \cite{G96},
\cite{G09, Ch.\ 10--11}), results on boundary value problems can be
adequately formulated within the scales of $H^s_q$- and
$B^s_q$-spaces (or just $H^s$-spaces) over $\Omega $ and $\partial\Omega $. For fractional-order pseudodifferential operators --- where the
prominent  example is the fractional Laplacian $(-\Delta )^a$,
$0<a<1$ --- we also need to introduce the $\mu $-transmission spaces
$H_p^{\mu (s)}(\comega)$, since they are the exact solution spaces for
Dirichlet problems.  The definition of these spaces involves a special type of
$\psi $do's called order-reducing operators.

The simplest examples of such operators (relevant for $\Omega =\rnp$)
are, with $t\in\R$,
$$
\Xi _\pm^t =\operatorname{Op}(\chi _\pm^t),\quad \chi _\pm^t(\xi )=(\ang{\xi '}\pm i\xi _n)^t ;\tag2.15
$$
they preserve support
in $\crnpm$, respectively, because the symbols extend as holomorphic
functions of $\xi _n$ into ${\Bbb C}_\mp$, respectively; here ${\Bbb
C}_\pm=\{z\in{\Bbb C}: \operatorname{Im}z\gtrless 0\}$. (The functions
$(\ang{\xi '}\pm i\xi _n)^t $  satisfy only part of the estimates
(2.9) (with $m=t$, $\tau \in\N_0$), but the $\psi $do definition can be applied anyway.)
 There is a more refined choice $\Lambda _\pm^t $, cf.\
\cite{G90, G15a}, with
symbols $\lambda _\pm^t (\xi )$ that do
satisfy all the required $\psi $do estimates, and where $\overline{\lambda _+^t
}=\lambda _-^{t }$.
These symbols likewise have holomorphic extensions in $\xi _n$ to the complex
halfspaces ${\Bbb C}_{\mp}$, so that  
the operators preserve
support in $\crnpm$, respectively. Operators with that property are
called ``plus'' resp.\ ``minus'' operators.
There is also a pseudodifferential definition $\Lambda
_\pm^{(t )}$ adapted to the situation of a bounded smooth domain $\Omega
$, by \cite{G90,G15a}.

The operators define homeomorphisms 
$
\Xi^t _\pm\colon H_q^s(\R^n) \simto H_q^{s- t
}(\R^n)$, for all $s\in \R$. 
The special
interest is that the ``plus''/``minus'' operators also 
 define
homeomorphisms related to $\crnp$ and $\comega$, for all $s\in{\Bbb R}$: 
$$
\aligned
\Xi ^{t }_+\colon \dot H_q^s(\crnp )&\simto
\dot H_q^{s- t }(\crnp),\quad
r^+\Xi ^{t }_{-}e^+\colon \ol H_q^s(\rnp )\simto
\ol H_q^{s- t } (\rnp ),\\
\Lambda  ^{(t) }_+\colon \dot H_q^s(\comega )&\simto
\dot H_q^{s- t }(\comega),\quad
r^+\Lambda  ^{(t) }_{-}e^+\colon \ol H_q^s(\Omega  )\simto
\ol H_q^{s- t } (\Omega  ),
\endaligned\tag2.16
$$
with similar rules for $\Lambda ^t_\pm$.
Moreover, the operators $\Xi ^t _{+}$ and $r^+\Xi ^{t }_{-}e^+$ identify with each other's adjoints
over $\crnp$, because of the support preserving properties.
There is a
similar statement for $\Lambda ^t_+$ and  $r^+\Lambda ^t_-e^+$
relative to $\rnp$, and for  $\Lambda ^{(t )}_+$ and $r^+\Lambda ^{(
t )}_{-}e^+$ relative to the set $\Omega $.  (The exponent $t$, and
the value $\mu $ considered below,  can also be allowed to be
complex, as in \cite{G15a,G22a}, but we shall not need this in the
present paper.) There is an abbreviation  in $\psi $do-notation $P_+=r^+Pe^+$ that is often
used, e.g.\ replacing $r^+\Xi _-^\mu e^+$ by $\Xi ^\mu _{-,+}$.

Let $\mu >-1$. Then the $\mu $-transmission spaces $H_q^{\mu
(s)}(\crnp)$ are defined for $s>\mu -1/q'$ by
$$
H_q^{{\mu} (s)}(\crnp)=\Xi _+^{-{\mu} }e^+\ol H_q^{s- {\mu}
}(\rnp).\tag2.17
$$
(Equivalently, $\Xi _+^{-\mu }$ can be replaced by  $\Lambda  _+^{-\mu
}$.) For $\mu =a>0$, the interest is that this is {\it the solution space}
for the homogeneous Dirichlet problem 
$$
(1-\Delta )^au=f\text{ on }\rnp,\quad \supp u\subset \crnp,
$$
  when $f$ is given in $\ol H_q^{s-2a}(\rnp)$ for some $s>a-1/q'$ and
  $u$ is sought in $\dot H_q^{\sigma}(\crnp)$ with $\sigma >a-1/q'$.
By \cite{G15a}, the result holds also for suitable
  variable-coefficient $\psi $do generalizations of 
$(1-\Delta )^a$. A pedestrian introduction to transmission spaces is
  given in \cite{G22c}.

The symbol $(\ang{\xi '}+i\xi _n)^{-\mu }$ is connected
  to expressions with a factor $x_n^\mu $ on $\crnp$ by the formula
$$
{\Cal F}^{-1}_{\xi _n\to x_n}\frac1{(\ang{\xi '}+i\xi _n)^{\mu +1}}=\tfrac 1{\Gamma (\mu +1)}e^+r^+x_n^\mu e^{-\ang{\xi '} x_n},\tag2.18
$$
 which allows to show
$$
H_q^{\mu (s)} (\overline{\Bbb R}^n_+)
\cases
=\dot H_q^{s}(\overline{\Bbb R}^n_+)\text{ if }\mu -1/q'<s<\mu +1/q,\\
\subset\dot
H_q^{s\,(-\varepsilon )}(\overline{\Bbb R}^n_+)+ e^+ x_n^\mu \overline H_q^{s-\mu }({\Bbb R}^n_+)\text{ if }s>\mu +1/q\endcases\tag2.19
$$
(with $(-\varepsilon )$ active if $s-\mu -1/q\in\N$).

For smooth bounded domains $\Omega $, the $\mu $-transmission spaces  $H^{\mu (s)}(\comega)$ are
defined in a similar way as in (2.17) with the operator family
$\Lambda _+^{(t)}$ used instead of $\Xi _+^{t}$; then there is a similar
inclusion as in (2.19) with $x_n^\mu $ replaced by $d_0^\mu $, and
these spaces are the solution spaces for homogeneous Dirichlet
problems for a large class of $\psi $do's \cite{G15a, Th.\ 4.4}
(containing the case $(-\Delta )^a$, $\mu =a$). There
is an analysis describing the spaces with further precision in \cite{G19}.

For $C^{1+\tau
}$-domains $\Omega $, the $\mu $-transmission spaces  $H^{\mu (s)}(\comega)$ are defined in
\cite{AG21, Def. 4.2}, when $\tau >0$, $\mu >-1$ and $\mu
-1/q'<s<1+\tau $, by localization, i.e., a reduction to local
coordinates in a family of open sets covering the boundary. When $\tau
\ge 1$, one then has with $\varepsilon >0$ \cite{AG21,Th.\ 4.5}:
$$
H_q^{{\mu} (s)}(\comega)\cases   =\dot H_q^{s}(\comega ),&\text{when }  s< {\mu} +\tfrac1{q},\\
 \subset \dot H_q^{s-\varepsilon }(\comega ), &\text{when }  s= {\mu} +\tfrac1{q},\\
 \subset  \dot H_q^{s }(\comega)+ d_0^\mu e^+\ol
  H_q^{s-\mu}({\Omega}),& \text{when }
  s-{\mu} -\tfrac1{q}\in \rp\setminus \N,\\
   \subset  \dot H_q^{s-\varepsilon }(\comega)+ d_0^\mu e^+\ol
  H_q^{s-\mu}({\Omega})& \text{when }
  s- {\mu} -\tfrac1{q}\in\N;
\endcases\tag2.20
$$
it also holds
 with $d_0$ replaced by $d$ (cf.\ (2.2)). When $\tau <1$,
there is a version of (2.20) with $d_0$ replaced by  local
choices of $d$, cf.\
\cite{AG21, Rem.\ 4.6}; for convenience, we recall the explanation
 here:

\example{Remark 2.1} When $\Omega $ is bounded  $C^{1+\tau}$-domain, each point $x_0\in
\partial\Omega $
has a bounded open neighborhood $U\subset \R^n$ and
a function $\zeta  \in C^{1+\tau }(\R^{n-1})$,
such that (after a suitable rotation) $\Omega \cap U=\R^n_\zeta  \cap
U$ (cf.\ (2.1)). For $\R^n_\zeta $, the space $H^{\mu(s)}_q(\ol{\R}^n_\zeta )$ is
defined from $ H^{\mu(s)}_q(\crnp)$ by use of the $C^{1+\tau
}$-diffeomorphism $F_\zeta (x)=
(x',x_n-\zeta (x'))$ (all $x\in\R^n$), with the notation $ F_\zeta
^{\ast,-1} ( u)\equiv   u\circ F_\zeta  ^{-1}$. Then
$H^{\mu(s)}_q(\ol{\R}^n_\zeta  ) = F^*_\zeta  (H^{\mu(s)}_q(\crnp))$,
provided with the inherited norm.

Now
$H^{\mu(s)}_q(\overline{\Omega})$ is defined as the set of all $u\in
H^s_{q,loc}(\Omega)$ such that for each $x_0\in\partial\Omega $, with a $\varphi\in
C_0^\infty(U)$ with $\varphi\equiv 1$ in a neighborhood of $x_0$, we have
$
  F_\zeta  ^{\ast,-1} (\varphi u)
  \in H^{\mu(s)}_q(\ol{\R}^n_+)$,
in the rotated situation.
The analysis of properties of
$H^{\mu(s)}_q(\overline{\Omega})$ is then carried over to an analysis of properties of
$H^{\mu(s)}_q(\overline\R^n_\zeta )$. The choices of $\{x_0,U,\zeta
,\varphi \}$ can be reduced to a finite system $\{x_{0,i},U_i,\zeta
_i,\varphi _i\}_{i=1,\dots,I}$, where  $\bigcup U_i$ covers
$\partial\Omega $.  The definition of the spaces
$H^{\mu(s)}_q(\overline\R^n_{\zeta_i} )$ is then used in each of the
sets $U_i$; and they may be  pieced together by a partition of unity. One has a
version of (2.20) with $\comega $ replaced by $\ol{\R}^n_{\zeta _i}$,
$d_0$ replaced by $d_i(x)=x_n-\zeta _i(x')$ for each $i$. When $\tau
\ge 1$, $d_i$ can be replaced by the $C^{1+\tau }$-function $d_0$ for
each $i$, and the pieces sum up to give the function $d_0$ entering in
(2.20). When $\tau \in \,]0,1[\,$, we still have the version of (2.20)
with $d_i$ in each localized piece $U_i$, where $d_i$ is related to
$d_0$ by (2.2), but $d_0$ is only $C^\tau $. 
\endexample

For $s\ge \mu +\frac1q$ there is also a weaker result than the second
line in (2.19), namely,
$$
H_q^{\mu (s)}(\crnp)=\Xi _+^{-\mu }e^+\ol H_q^{s-\mu }(\rnp )\subset
\Xi  _+^{-\mu }\dot H_q^{\frac1q-\varepsilon }(\crnp)=\dot H_q^{\mu
+\frac1q-\varepsilon }(\crnp), \text{ when }s\ge \mu +\tfrac1q.
$$
In the case of $C^{1+\tau }$-domains, the corresponding rule
 $$
H_q^{\mu (s)}(\comega)\subset \dot H_q^{\mu +\frac1q-\varepsilon }(\comega
),\text{ for }s\ge \mu +\tfrac1q,\;  s<1+\tau  ,\tag 2.21
$$
follows by use of 
 the localization described in Remark 2.1.

\subsubhead 2.4 Trace mappings, additional properties \endsubsubhead

The weighted trace mapping
$$
 \gamma
  _0^\mu \colon u\mapsto \Gamma (\mu +1)(u/d^\mu )|_{\partial \Omega  },\tag2.22
$$
is defined on $H_q^{\mu (s)}(\comega)$ for $C^{1+\tau }$-domains in \cite{AG21, Sect.\ 4.2}. 
Since we shall in the present paper deal with nonhomogenous boundary
values, we need a more elaborate version of some results from there.
The statement on $\gamma _0^\mu $ in \cite{AG21, Prop.\
4.3} for a curved halfspace (2.1) extends as follows:

\proclaim{Proposition 2.2} Let $\mu >-1$, $\tau >0$, and $\mu +\frac1{q}< s<1+\tau $ with
  $s-\mu<1+\tau$, 
  and let $\rn_\zeta $ be defined by  $\zeta \in C^{1+\tau}(\R^{n-1})$, $d(x)=x_n-\zeta (x')$ near
  $\partial\rn_\zeta $.  The mapping $\gamma
  _0^\mu \colon u\mapsto \Gamma (\mu +1)(u/d^\mu )|_{\partial
           \rn_\zeta }$
  is continuous and surjective:
$$
\gamma _0^\mu \colon H_q^{{\mu} (s)}(\ol{\R}^n_\zeta )\to B_{q}^{s-\mu
  -\frac1q}(\partial \rn_\zeta),\tag2.23
$$
having a continuous right inverse. Moreover, the space $H_q^{{(\mu+1)}
(s)}(\ol{\R}^n_\zeta  )$ is a closed subspace of $H_q^{{\mu} (s)}(\ol{\R}^n_\zeta  )$,
equal to the kernel of the mapping {\rm (2.23)}.
\endproclaim

\demo{Proof} The properties are known from \cite{G15a, Sect.\ 5} to hold when $\zeta (x')\equiv 0$
(the flat case), as recalled e.g.\ in \cite{AG21, Sect.\ 4.1, (4.7)}. They carry over to
the case of general $\zeta (x')$ in view of the  mapping properties of the
diffeomorphism $F_\zeta \colon (x',x_n)\mapsto (x',x_n-\zeta (x'))$ and the
definitions 
listed in \cite{AG21, Sect.\ 4.2}. Note that the space $H_q^{{(\mu+1)}
(s)}(\ol{\R}^n_\zeta  )$ is well-defined with the parameter $\mu +1$, since
the hypotheses assure that  $s>(\mu +1)-1/{q'}$.\qed  
\enddemo

When $\tau \geq 1$  one can, as stated in \cite{AG21, Prop.\ 4.3},
 replace $d$ by $d_0$.
We then also have the following elaborated version of the statements
on $\gamma _0^\mu $ for bounded domains $\Omega $ in
\cite{AG21, Th.\ 4.5}: 

\proclaim{Theorem 2.3} Let $\mu >-1$, $\tau \ge 1$ and $\mu +\frac1{q}< s<\tau $ with
  $s-\mu<\tau$, and  let $\Omega\subset\R^n$ be a bounded
  $C^{1+\tau}$-domain.
  The mapping $\gamma _0^\mu \colon u\mapsto \Gamma (\mu +1)(u/d_0^\mu
)|_{\partial \Omega }$ is continuous and surjective:
$$
\gamma _0^\mu \colon H_q^{{\mu} (s)}(\comega )\to B_q^{s-\mu
  -\frac1q}(\partial \Omega ),\tag2.24
$$
having a  continuous right inverse $K^{\mu }_{(0)}$. Moreover, the space $H_q^{{(\mu +1)}
(s)}(\comega )$ is a closed subspace of $H_q^{{\mu} (s)}(\comega )$,
and equals the kernel of the mapping {\rm (2.24)};
$$
 \{u\in H_q^{\mu (s)}(\comega)\mid\gamma
_0^{\mu }u=0\}=H_q^{(\mu +1)(s)}(\comega). \tag2.25
$$

\endproclaim

\demo{Proof} The continuity of the mapping $\gamma _0^\mu $ is
established in \cite{AG21, Th.\ 4.5} by use of a cover
$\bigcup_{i=0,1,\dots, I}U'_1$ and an associated partition of unity
$\{\varrho _i\}_{i=0,\dots,I}$ such that the $U'_i$ with $i\ge 1$
cover $\partial\Omega $ and for each such $U'_i$ there is a function
$\zeta _i$ (called $\gamma _i$ in \cite{AG21}) such that, after a
rotation and translation depending on $i$, $\Omega \cap U'_i =
\rn_{\zeta _i}\cap U'_i$. In each such neighborhood $U'_i$, the facts
known for $\rn_{\zeta _i}$ can be applied to $\varrho _iu$ and
collected to a statement on $u$ by summation. This goes for all the
properties listed in Proposition 2.2, when we moreover observe that
$\gamma _0^\mu $ acts locally as a trace operator ($\gamma
_0^a(\varphi u)=\gamma _0\varphi \gamma _0^au$ when $\varphi \in
C_0^\infty (\rn)$), and $d_0$ is defined
near the boundary consistently in the different local charts. \qed
\enddemo

\remark{Remark 2.4} For a description of $K^\mu _{(0)}$, we note that the analysis in \cite{G19,
Sect.\ 3} for smooth $\Omega $ shows that in local coordinates reduced to the case of
$\rnp$, $K_{(0)}^{\mu }$ can be taken proportional to 
$x_n^{\mu }K_0$, where $K_0$ is the standard Poisson operator for the
Laplacian (solving the problem $(1-\Delta )
v=0$ on $\rnp$, $\gamma _0v=\varphi $). This explains the factor $d^\mu $
appearing when the result is carried over to $\Omega $.
\endremark

One can also define higher-order traces $\gamma _k^\mu u$;
$$
\gamma _k^\mu u=\Gamma (\mu +k+1) \gamma _k(u/d^\mu ),\tag 2.26
$$
when $\mu +\frac1q+k<s<\tau $. Since they are not used in this paper,
we leave out details.

Let us finally recall the scale of spaces built over the H\"older-Zygmund
spaces $C^s_*(\rn)$, with completely parallel properties, as accounted
for in \cite{AG21}. With the associated scales over $\ol C_*^s(\rnp)$
and $\dot C_*^s(\crnp )$ defined
by (2.6), one defines $C_*^{\mu (s)}(\crnp)$ (when $\mu >-1$) similarly to (2.17) by
$$
C_*^{{\mu} (s)}(\crnp)=\Xi _+^{-{\mu} }e^+\ol C_*^{s- {\mu}
}(\rnp), \quad s>\mu -1.\tag2.27
$$
For bounded $C^{1+\tau }$-domains, corresponding spaces are defined by localization
when $\tau >0$, $\mu -1<s<1+\tau $, and have the properties when $\tau
\ge 1$:
$$
C_*^{{\mu} (s)}(\comega)\cases   =\dot C_*^{s}(\comega ),&\text{for }  s< {\mu} ,\\
 \subset  \dot C_*^{s\,(-\varepsilon )}(\comega)+ d_0^\mu e^+\ol
  C_*^{s-\mu}({\Omega}),& \text{for }
  s> {\mu} 
\endcases\tag2.28
$$
(with $(-\varepsilon )$ active if $s-\mu \in\N$), cf.\ \cite{AG21,
Def.\ 4.2, Th,\ 4.5}. This also holds
with $d_0$ replaced by $d$ (cf.\ (2.2)). When $\tau <1$,
there is a version of (2.28) with $d_0$ replaced by  local
choices of $d$, cf.\
 Remark 2.1. We also have:
$$
C_*^{\mu (s)}(\comega)\subset \dot C_*^{\mu -\varepsilon }(\comega
),\text{ for }s> \mu >0,\;  s<1+\tau  .\tag 2.29
$$

The trace mapping   $\gamma _0^\mu \colon u\mapsto \Gamma (\mu +1)(u/d_0^\mu
)|_{\partial \Omega }$
$$
\gamma _0^{\mu }\colon C_*^{{\mu} (s)}(\comega )\to C_*^{s-\mu
  }(\partial \Omega ),\tag2.30
$$
is well-defined when 
$s>\mu $. Since it plays an important
role in the study of nonhomogeneous boundary problems, we include an
elaborated version of the statements on $\gamma _0^\mu $ in
\cite{AG21}, along the lines of Theorem 2.3:

\proclaim{Theorem 2.5} Let $\mu >-1$, $\tau \ge 1$ and $\mu < s<\tau $ with
  $s-\mu<\tau$, and  let $\Omega\subset\R^n$ be a bounded
  $C^{1+\tau}$-domain.
  The mapping $\gamma _0^\mu \colon u\mapsto \Gamma (\mu +1)(u/d_0^\mu
)|_{\partial \Omega }$ in {\rm (2.30)} is continuous and surjective,
having a  continuous right inverse. Moreover, the space $C_*^{{(\mu +1)}
(s)}(\comega )$ is a closed subspace of $C_*^{{\mu} (s)}(\comega )$,
and equals the kernel of the mapping {\rm (2.30)}.
\endproclaim

The proof goes exactly as in the proof of Theorem 2.3.

As noted in \cite{AG21, Cor.\ 6.11}, the well-known embedding property $H_q^t(\rn)\subset
C_*^{t-n/q-\varepsilon }(\rn)$ (any $\varepsilon >0$) implies
$$
H_q^{{\mu} (s)}(\crnp)\subset C_*^{{\mu} (s-n/q-\varepsilon )}(\crnp)
$$
in view of the
definitions (2.17) and (2.27); this leads to  embeddings
$$
H_q^{{\mu} (s)}(\comega)\subset C_*^{{\mu} (s-n/q-\varepsilon )}(\comega),\tag2.31
$$
when $\mu -1 < s-n/q-\varepsilon <s<1+\tau $ and $\Omega $ is a bounded
$C^{1+\tau }$-domain, by following the localization procedure.
Letting $q\to\infty $, we find for $\mu -1<s<1+\tau $:
$$
\bigcap_{q>1}H_q^{{\mu} (s)}(\comega)\subset C_*^{{\mu} (s-\varepsilon )}(\comega),\tag2.32
$$
for small $\varepsilon >0$.
In the other direction, one has e.g.\ for $0<t<1+\tau$, 
$$
C_*^{t+\varepsilon }(\comega)\subset \ol H_q^{t}(\Omega ).\tag2.33
$$
These observations are useful for drawing consequences on regularity in
H\"older spaces.

\subhead 3. Mapping properties and regularity \endsubhead

\subsubhead 3.1 The homogeneous Dirichlet problem
\endsubsubhead

In the rest of this paper, we consider a classical pseudodifferential
 operator
 $P$ of order $2a$, $0<a<1$, with {\it
even} symbol in
$C^\tau S^{2a}(\rn\times\rn)$, cf.\ (2.14). As accounted for in
 \cite{AG21}, the so-called $a$-transmission condition for an operator
 of order $2a$, with respect to
 a $C^{1+\tau }$-domain $\Omega $, means that
 (2.14) is satisfied at each $x\in\partial\Omega $ for $\xi $ equal
 to the interior normal $\nu (x)$ at $x$. The evenness implies that $P$ satisfies
the $a$-transmission condition with respect to any domain, and it is
 shown how the spaces
$H_q^{a(t)}(\comega)$ enter in discussions of mapping properties. We remark that $P$ also satisfies the $(a+k)$-transmission
 condition for $k\in 
\Z$, since a $\mu $-transmission condition only depends on $\mu $
modulo 1.

Recall that 
we assume
$1<q<\infty $ throughout. Our assumptions for the treatment of the homogeneous
Dirichlet problem are as follows:

\proclaim{Hypothesis 3.1}
 $1^\circ$ There are given constants $a,\tau ,q$ with $0<a<1$, $\tau >2a$, and 
 $1<q<\infty $.  
 $\Omega $ is a bounded $C^{1+\tau }$-domain in $\rn$, and  $P$ is a
 classical $\psi $do of order $2a$, with  even symbol in
$C^\tau S^{2a}(\rn\times\rn)$.

$2^\circ$ Assumptions as in $1^\circ$, and in addition $P$ is strongly elliptic.

\endproclaim

The following results for the homogeneous Dirichlet problem were shown
in \cite{AG21} (Theorems 6.4 and 6.9, Cor.\ 6.10):

\proclaim{Theorem 3.2} Assume Hypothesis {\rm 3.1 $1^\circ$}, and let
$s$ satisfy  $-a\le s<\tau -2a$. Then
$r^+P$ maps continuously
$$
   r^+ P\colon H_q^{a (s+2a)}(\comega)\to \ol H_q^{s}(\Omega ).\tag3.1
   $$

 Assume moreover that Hypothesis {\rm 3.1 $2^\circ$} holds and  $s\ge 0$. If
 $u\in
 \dot{H}^a_q(\ol{\Omega})$ solves
 $$
     Pu = f \text { in }\Omega ,\quad \supp u\subset \comega,\tag3.2
 $$
  for some $f\in \ol{H}^s_q(\Omega)$, then $u\in
  H^{a(s+2a)}_q(\ol\Omega)$.

\endproclaim

In other words, we have found that
$$
\{u\in \dot H_q^a(\comega)\mid r^+Pu\in \ol H_q^s(\Omega )\}=H^{a(s+2a)}_q(\comega).\tag3.3
$$

 Note the remarkable fact that we have not only shown a regularity
    property (a conclusion from $Pu$ to $u$) as the aim is  in
    most studies of these operators, but we have found {\it the exact solution space}
    for the homogeneous Dirichlet problem with data in $\ol
    H_q^s(\Omega )$; the {\it Dirichlet domain}. Moreover,  $H^{a(s+2a)}_q(\ol\Omega)$ has this
    role {\it universally}, in the sense that
it is independent of the choice of $P$ satisfying Hypothesis 3.1.

Let us derive some consequences of the theorem, particularly concerning
the range of  $r^+P$:

\proclaim{Proposition 3.3} Assume Hypothesis {\rm 3.1} and
let $0\le s<\tau -2a$. The solutions of {\rm (3.2)} described in
Theorem {\rm 3.3} satisfy an estimate
$$
\|u\|_{H^{a(s+2a)}_q(\ol\Omega)}\le
C(\|r^+Pu\|_{\ol{H}^s_q(\Omega)}+\|u\|_{
\dot{H}^a_q(\ol{\Omega})}). \tag 3.4
$$
The range of $r^+P$ in {\rm (3.1)} is closed. 
\endproclaim

\demo{Proof}
The estimate (3.4) is seen as follows: 
By the continuity statement in Theorem 3.2, $|\|u|\|\equiv \|r^+Pu\|_{\ol{H}^s_q(\Omega)}+\|u\|_{
\dot{H}^a_q(\ol{\Omega})}$ is a norm on $H_q^{a(s+2a)}(\comega)$ that
satisfies $|\|u|\|\le C'\|u\|_{H^{a(s+2a)}_q(\ol\Omega)}$. By use of the
regularity statement, ${H^{a(s+2a)}_q(\ol\Omega)}$ is seen to be complete under
this norm. Then since it is a Banach space, the norms are equivalent.

The closed range property follows from a classical argument for how an
a priori
inequality (3.4) leads to the existence of an approximate inverse,
found e.g.\ in the proof of \cite{H63, Th.\ 10.5.1}. Denote
$\ol H_q^s(\Omega )=X$, $ H_q^{a(s+2a)}(\comega)=Y$, $ \dot H_q^a(\comega)=Z$,
they are all reflexive Banach spaces. Here $Y$ is compactly injected in
$Z$, since $ H_q^{a(s+2a)}(\comega)\subset \dot
H_q^{a+b}(\comega)\subset \dot H_q^a(\comega)$ for any
$0<b<\min\{a,\tfrac1q\}$ (cf.\ (2.20)--(2.21)), where the last injection is compact. Write
$r^+P$ as $P$ for short. For $u\in
Y$, 
the inequality (3.4) reads
$$
\|u\|_Y \le  C(\|Pu\|_X + \|u\|_Z). \tag 3.5
$$
Let $N$ denote the nullspace of $P\colon Y\to X$; it is a closed subspace of $Y$,
and there is a closed complement $Y'$ of $N$ in $Y$ such that $P\colon
Y'\to X$ is injective (and bounded). Let $R$ denote the range of $P$
in $X$. Our goal will be achieved if we show that
$$
\|u\|_{Y'} \le  C'\|Pu\|_X \text{ when }u\in Y'; \tag 3.6
$$
 for then $P\colon Y'\to X$ has a bounded partial inverse $Q\colon
 R\to Y'$; and since $Q$ is closed, $R$ is a closed subspace of $X$.

Assume that (3.6) does not hold. Then there is a sequence of
functions $u_k\in Y'$ with $\|u_k\|_Y=1$ such that $\|Pu_k\|_X \to
0$. By the weak compactness of the unit sphere in $Y$, there is a
subsequence converging weakly to an element $u_0$. By the compactness
of the injection $Y \subset Z$, we can take a further subsequence (call it $u_k$ again) such that $\|u_k-u_0\|_Z \to 0$. Now (3.5) implies
$$
\|u_k-u_0\|_Y \le  C(\|P(u_k-u_0)\|_X + \|u_k-u_0\|_Z).\tag3.7
$$
Since $P$ is a fortiori weakly continuous from $Y$ to $X$, $Pu_k$ goes
weakly to $ Pu_0$, which must be $0$ since $Pu_k$ goes to $0$ in
norm. Then both terms in the right-hand side of (3.7) go to $0$. If $u_0=0$, this
gives a contradiction since $u_k$ has norm 1. If $u_0\ne 0$, we have
found a nontrivial null-element in $Y'$, contradicting the definition of
$Y'$. Thus  (3.6) must hold.\qed

\enddemo

We shall later (in Section 4) show a bijectiveness or Fredholm property of the mapping
(3.1), leading to existence and uniqueness results for the homogeneous
Dirichlet problem. It is important for that study to observe that not
only the regularity parameter $s$, but also the integral
parameter $q$ for $u$ can be lifted when the data are in a space with a higher 
parameter:

\proclaim{Theorem 3.4} Assume Hypothesis {\rm 3.1}.

If $u\in \dot H_q^a(\comega)$ solves {\rm (3.2)} for some $f\in
L_p(\Omega )$ with $p\ge q$, then $u\in H_p^{a(2a)}(\comega)$.

In fact, there is a sequence $q_1=q<q_2<q_3<\cdots$ with $q_j\to\infty
$ for $j\to\infty $, such that for all $j$, $u\in \dot
H^a_{q_{j}}(\comega)$ with  $r^+Pu\in L_{q_{j+1}}(\Omega )$ imply $u\in
\dot H_{q_{j+1}}^a(\comega)$.
\endproclaim

\demo{Proof}
The integral parameter $q$ will be lifted by use of the well-known embedding rule
$$
H^{s_1}_{p_1}(\rn)\subset H^{s_2}_{p_2}(\rn)\text{ when }s_1\ge s_2,\;
s_1-\tfrac n{p_1} \ge s_2-\tfrac n{p_2}, \; 1< p_1\le p_2<\infty . \tag3.8
$$
It holds also for $\dot H^s_p(\comega )$-spaces with $\Omega$ open
$\subset \rn$, by their definition as closed subspaces of
$H^s_p(\rn)$, cf.\ (2.5).
 
In
view of (2.20) and (2.21) we have for $1<p<\infty $,
 $$
 H_p^{a(2a)}(\comega)\cases = \dot
H_p^{2a}(\comega)\text{ if }a<\frac1p\\
\subset\dot H_p^{a+\frac1p-\varepsilon }(\comega)\text{ if
}a\ge \tfrac1p.
\endcases\tag3.9
$$
To reduce this to one statement, define
$$
m(p)=\min \{ap,1\}\text{ for  }p\ge q;
$$
it is 1 when $p\ge 1/a$, and if $a<\frac1q$ so that the interval
between $q$ and $1/a$ is nontrivial, $m(p)$ takes values in $[aq,1]$
on that interval. Altogether, $m(p)$ satisfies
$$
\min\{aq,1\}\le m(p)\le 1\text{ for all }p\ge q.\tag3.10
$$
Then (3.9) implies
$$
H_p^{a(2a)}(\comega)\subset \dot H_p^{a+\frac1p m(p)-\varepsilon
}(\comega)\text{ when }p\ge q. 
$$
Now we define an increasing sequence  $q_1<q_2<\dots <q_j<\dots$,
where $q_1=q$, and the next values are defined successively such that
$$
 \dot H_{q_j}^{a+\frac1q m(q_j)-\varepsilon }(\comega)\subset \dot H_{q_{j+1}}^a(\comega). 
$$
 By (3.8), this is satisfied if $a+\frac1{q_j} m(q_j)-\varepsilon -\frac n{q_j}\ge
 a-\frac n{q_{j+1}}$, which may be rewritten as $\frac
 n{q_{j+1}}>\frac{n-m(q_j)}{q_j}$; i.e.,
$$
q_{j+1}<\tfrac n
{n-m(q_j)}q_j.\tag 3.11
$$
In view of (3.10),
the factor $n/(n-m(p))$ is for all $p\ge q$ greater than a fixed constant $c>1$, so the
inequality (3.11) allows defining $q_j$ as a sequence going to $\infty
$ (it holds e.g.\ if we take $q_j=c_1^{j-1}q_1$ with a $c_1\in \,]1,c[\,$, for $j\in\N$). Note that we have obtained
$$
H_{q_j}^{a(2a)}(\comega)\subset  \dot H_{{q_{j+1}}}^a(\comega),\text{
for all }j\in\N.\tag3.12
$$
The sequence satisfies the assertion in the theorem.

For the given $u\in \dot H_q^a(\comega)$ and $f\in L_p(\Omega )$, we
use this in a
successive passage from $q_j$ to $q_{j+1}$ for $j=1,2,\dots$, until
$p$ is reached: 
Starting with the
information $r^+Pu=f\in L_p(\Omega )\subset L_q(\Omega )$, $q_1=q$, we have from Theorem 3.2  that $u\in
H_{q_1}^{a(2a)}(\comega)$, hence is in  $ \dot H_{q_2}^{a}(\comega)$ by (3.12). If
$q_2\ge p$, this ends the proof since $u\in \dot H_{p}^{a}(\comega)$
then, so that by Theorem 3.2, $u\in
H_{p}^{a(2a)}(\comega) $. Otherwise, we repeat the argument. The
general step, as long as $q_j<p$, is that  the
information $r^+Pu=f\in L_p(\Omega )\subset L_{q_j}(\Omega )$ gives by Theorem 3.2  that $u\in
H_{q_j}^{a(2a)}(\comega)$, which is in $ \dot
H_{q_{j+1}}^{a}(\comega)$ by (3.12). If
$q_{j+1}\le p$,  Theorem 3.2 implies that $u\in
H_{q_{j+1}}^{a(2a)}(\comega)$. When $q_{j+1}>p$, Theorem 3.2 is applied
with the parameter $p$. \qed
\enddemo

\subsubhead 3.2 The nonhomogeneous local Dirichlet problem
\endsubsubhead

We shall now include  nontrivial Dirichlet boundary values, and therefore go out in the larger
spaces $H^{(a-1)(s+2a)}_q(\comega)$ to define solutions. These spaces are
defined as in Section 2.3 with $\mu =a-1$ (which is $>-1$).
To the basic assumptions listed in Hypothesis 3.1, we add the condition
$\tau >2a+1$.

There is
the following forward mapping property based on \cite{AG21, Cor.\ 5.14}:

\proclaim{Theorem 3.5} Asssume Hypothesis {\rm 3.1} $1^\circ$ with $\tau >2a+1$, and let
$s$ satisfy $-a-1\le s<\tau -2a-1$.
Then for curved halfspaces $\rn_\zeta $
with $\zeta \in C^{1+\tau }(\R^{n-1})$, and for bounded $C^{1+\tau
}$-domains $\Omega $,
$r^+P$ maps continuously
$$
\aligned
   r^+ P&\colon H_q^{(a-1) (s+2a)}(\ol{\R}^n_\zeta )\to \ol H_q^{s}({\R}^n_\zeta  ),\\
   r^+ P&\colon H_q^{(a-1) (s+2a)}(\comega)\to \ol H_q^{s}(\Omega ).
\endaligned
\tag3.13
$$
\endproclaim

\demo{Proof} By \cite{AG21, Cor.\ 5.14}, the first statement holds in
the flat case $\zeta \equiv 0$, since $P$ satisfies the global
$(a-1)$-transmission condition. It carries over to general $\zeta $
exactly as in the proofs of \cite{AG21, Th.\ 6.1, Cor.\ 6.2}, with $a$
replaced by $a-1$ in the definitions of transmission spaces.

For the second statement, this is carried over to bounded $C^{1+\tau
}$-domains analogously to the proof of \cite{AG21, Th.\ 6.4}. \qed
\enddemo
 Combining this with the mapping properties of $\gamma _0^{a-1}$ shown
 in Proposition 2.2 and Theorem 2.3, we have as a corollary:

\proclaim{Corollary 3.6} Assume  Hypothesis {\rm 3.1} $1^\circ$ with $\tau >2a+1$, and let
$s$ satisfy $-a-1/q'< s<\tau -2a-1$. 
The following maps are continuous:
$$
\aligned
   \{r^+ P,\gamma _0^{a-1}\}&\colon H_q^{(a-1) (s+2a)}(\ol{\R}^n_\zeta
   )\to \ol H_q^{s}({\R}^n_\zeta  )\times B_q^{s+a+1/q'}(\partial\rn_\zeta ),\\
   \{r^+ P,\gamma _0^{a-1}\}&\colon H_q^{(a-1) (s+2a)}(\comega
   )\to \ol H_q^{s}(\Omega   )\times B_q^{s+a+1/q'}(\partial\Omega  ),
\endaligned
\tag3.14
$$

\endproclaim

These results allow a development of regularity results for the nonhomogeneous Dirichlet
problem:
 $$
 \aligned
    Pu& = f \text { in }\Omega ,\\
     u&=0 \text { in }\rn\setminus \Omega ,\\
     \gamma _0^{a-1}u&=\varphi \text{ on }\partial\Omega .
     \endaligned \tag3.15
 $$

\proclaim{Theorem 3.7} Assume  Hypothesis {\rm 3.1}  with $\tau >2a+1$, and let
$s$ satisfy $0\le  s<\tau -2a-1$.

When $u\in H_q^{(a-1)(a)}(\comega)$,  it satisfies {\rm (3.15)}
 for some $f\in \ol H_q^{-a}(\Omega)$ and $\varphi \in
 B_q^{1/q'}(\partial\Omega  ) $.
 
Let $u\in H_q^{(a-1)(a)}(\comega)$, and let $p\ge q$. If $u$ solves {\rm (3.15)} with  $f\in \ol{H}^s_p(\Omega)$ and $\varphi
\in B_p^{s+a+1/p'}(\partial\Omega  ) $, then $u\in
H^{(a-1)(s+2a)}_p(\ol\Omega)$.

In other words,
$$
\{u\in  H_q^{(a-1)(a)}(\comega)\mid r^+Pu\in \ol H_p^s(\Omega ),
\gamma _0^{a-1}u\in B_p^{s+a+1/p'}(\partial\Omega  )   \}=H^{(a-1)(s+2a)}_p(\comega).\tag3.16
$$

\endproclaim

\demo{Proof} The first statement on $u$ follows from the
mapping property in Corollary 3.6 for $s=-a$.

Now for the improvement of regularity and integral parameters: Assume $f\in \ol{H}^s_p(\Omega)$ and $\varphi
\in B_p^{s+a+1/p'}(\partial\Omega  ) $. By the surjectiveness of $\gamma _0^{a-1}$ in Theorem
2.2, there is a $v\in H_p^{(a-1)(s+2a)}(\comega)$ with $\gamma
_0^{a-1}v=\varphi $. By Theorem 3.5, $g=r^+Pv\in \ol{H}^s_p(\Omega)
$. Thus $w=u-v\in H_q^{(a-1)(a)}(\comega)$ and solves
 $$
 \aligned
     Pw& = f -g\text { in }\Omega ,\\
     w&=0 \text { in }\rn\setminus \Omega ,\\
     \gamma _0^{a-1}w&=0 \text{ on }\partial\Omega ,
     \endaligned \tag3.17
 $$
 where $f-g\in \ol{H}^s_p(\Omega)$.
 This is a homogeneous Dirichlet problem, and since $w$ is in the
 kernel of $\gamma _0^{a-1}$, (2.25) shows that $w\in
 H_q^{a(a)}(\comega)= \dot H_q^a(\comega)$. Then Theorems 3.2 and 3.4
 apply to show that since $f-g\in
 \ol{H}^s_p(\Omega)$, $w$ is in $  H_p^{a(s+2a)}(\comega)$, and hence
 in the larger space  $ H_p^{(a-1)(s+2a)}(\comega) $.
 Finally, $u=v+w\in  H_p^{(a-1)(s+2a)}(\comega)$.\qed
\enddemo

The existence and uniqueness of solutions will be taken up in Section 5.

\example{Remark 3.8}  It is possible to describe the spaces
$H_q^{a(s)}(\comega)$ and $H_q^{(a-1)(s)}(\comega)$ more precisely.
We already know that $H_q^{a(s)}(\comega)=\dot H_q^s(\comega)$ when
$s<a+\frac1q$, and $H_q^{(a-1)(s)}(\comega)=\dot H _q^s(\comega)$ when
$s<a-1+\frac1q=a-\frac1{q'}$, by (2.20) for $\mu =a$ resp.\ $a-1$.

For higher $s$ (up to $\min\{\tau ,\tau +\mu \}$), let $K_{(0)}^{\mu }$ be  a right inverse
of $\gamma _0^{\mu }$ as in Theorem 2.3. The last statement in Theorem 2.3 implies
$$
H_q^{\mu (s)}(\comega)= H_q^{(\mu +1)(s)}(\comega)\,\dot + \,K_{(0)}^{\mu }B_q^{s-\mu -1/q}(\partial\Omega
), \text{ when }s> \mu +\tfrac1q .\tag3.18
$$
In particular, we have with $\mu =a$:
$$
H_q^{a(s)}(\comega)=\dot H_q^{s}(\comega)\,\dot+\,K_{(0)}^{a}B_q^{s-a-1/q}(\partial\Omega
)\text{ for }s-a\in \,]\tfrac1q,1+\tfrac1q[\,, \tag3.19
$$
and with $\mu =a-1$:
$$
H_q^{(a-1)(s)}(\comega)=\cases \dot H_q^{s}(\comega)\,\dot+\,K_{(0)}^{a-1}B_q^{s-a+1/q'}(\partial\Omega
)\text{ for }s-a\in\,]-\tfrac1{q'},\tfrac1q[\,,\\
\dot H_q^{s}(\comega)\,\dot+\,K_{(0)}^{a}B_q^{s-a-1/q}(\partial\Omega
)\,\dot +\,K_{(0)}^{a-1}B_q^{s-a+1/q'}(\partial\Omega
)\text{ for }s-a\in\,]\tfrac1{q},1+\tfrac1q[\,.\endcases
\tag3.20
$$
The operators $K_{(0)}^a$ and $K_{(0)}^{a-1}$ provide coefficients
$d^a$ resp.\ $d^{a-1}$, cf.\ Remark 2.4. In particular, when
$a<\frac1q$, $H_q^{a(2a)}(\comega)=\dot H_q^{2a}(\comega)$ and $H_q^{(a-1)(2a)}(\comega)=\dot H_q^{2a}(\comega)\,\dot+\,K_{(0)}^{a-1}B_q^{a+1/q'}(\partial\Omega
)$;  and when  $a>\frac1q$, $H_q^{a(2a)}(\comega)=\dot
H_q^{2a}(\comega)\,\dot+\,K_{(0)}^{a}B_q^{a-1/q}(\partial\Omega) $ and
$H_q^{(a-1)(2a)}(\comega)=\dot
H_q^{2a}(\comega)$\linebreak $\dot+\,K_{(0)}^{a}B_q^{a-1/q}(\partial\Omega)
\,\dot+\,K_{(0)}^{a-1}B_q^{a+1/q'}(\partial\Omega
)$.

\endexample

By playing on the possibility to take $p$ very large in Theorem 3.7,
we can draw a consclusion on problems with data in H\"older spaces:

\proclaim{Corollary 3.9} Assume  Hypothesis {\rm 3.1}  with $\tau >2a+1$, and let
$s$ satisfy $0\le  s<\tau -2a-1$.

Let $u\in H_q^{(a-1)(a)}(\comega)$ for some $1<q<\infty $. If $u$ solves {\rm (3.15)} with  $f\in C^s(\comega)$ and $\varphi
\in C^{s+a+1}(\partial\Omega  ) $, 
then $$
u\in
C_*^{(a-1)(s+2a-\varepsilon )}(\ol\Omega)\subset \dot
C^{s+2a-\varepsilon }(\comega)+d^{a-1}C^{s+a+1-\varepsilon }(\comega),
$$ for small $\varepsilon >0$ (with noninteger $s+2a-\varepsilon $ and $s+a+1-\varepsilon $).
\endproclaim

\demo{Proof} In view of (2.33), we have for 
large $p\in [q,\infty [\,$ and small $\varepsilon '>0$ 
that $f\in \ol H_p^{s-\varepsilon '}(\Omega )$ and $\varphi \in
B_p^{s+a+1/p'-\varepsilon '}(\partial\Omega )$. Then Theorem 3.7
implies that $u\in H_p^{(a-1)(s-\varepsilon '+2a)}(\comega)\subset \dot
H_p^{s-\varepsilon '+2a}(\comega)+d^{a-1}e^+\ol H_p^{s-\varepsilon
'+a+1}(\Omega )$. The result follows in view of (2.32) by letting
$p\to\infty $. \qed
\enddemo

\subhead 4. Eigenfunctions and Fredholm properties \endsubhead

We shall now derive existence-and-uniqueness (or Fredholm) properties of
$r^+P$ acting as in (3.1). This will follow from an analysis of
eigenfunctions, also for $r^+P^*$, of the Dirichlet realizations they define in
$L_q(\Omega )$.

\subsubhead 4.1 The Dirichlet realization in $L_2(\Omega )$
\endsubsubhead

In order to define the Dirichlet realization in $L_2(\Omega )$ of a
strongly elliptic pseudodifferential operator $P$, we  first recall the G\aa{}rding inequality, adapting it to the present symbols.

\proclaim{Lemma 4.1} When $\tau >a>0$, 
and $P=\Op(p)$ with symbol $p(x,\xi )$ in $ C^\tau
S^{2a}(\rn\times\rn)$ is strongly elliptic, i.e.\, the
principal symbol satisfies $$
\operatorname{Re}p_0(x,\xi )\ge c_0|\xi
|^{2a}\text{ for all }x\in \rn, |\xi |\ge 1,\tag4.1
$$ then the operator satisfies
the G\aa{}rding inequality:
$$
\operatorname{Re}(Pu,u)\ge c\|u\|^2_{
H^a(\rn)}-\beta \|u\|^2_{L_2(\rn )}, \text{ for }u\in C_0^\infty (\rn),\tag 4.2
$$
for some $ c>0$ and $\beta \in{\Bbb R}$. One can take $c=c_0-\delta $, any
$\delta \in \,]0,c_0[\,$.
\endproclaim

\demo{Proof} The G\aa{}rding inequality for smooth symbols is an old
and well-established fact. The reader who wants to see a proof on
$\rn$ can find it e.g.\ in \cite{G11, Lemma 5.1}, which carries over
verbatim to the case where $A$ is replaced by a $\psi $do $P$ with
symbol  $p$ in
$S^{2a}(\rn\times\rn)$ ($p_0$ taken smooth near $\xi =0$ with $\operatorname{Re}p_0>0$), and $m$ is replaced by $a$.

Note here that
since $P-(c_0-\delta )\Op(\ang\xi ^{2a})$ is likewise 
strongly elliptic when $\delta \in \,]0,c_0[\,$, 
$$
\aligned
\operatorname{Re}((P-(c_0-\delta )\Op(\ang\xi
^{2a}))u,u)&=\operatorname{Re}(Pu,u)-(c_0-\delta )\|u\|^2_{H^a(\rn)}\\
&\ge c'\|u\|^2_{
H^a(\rn)}-\beta '\|u\|^2_{L_2(\rn )},
\endaligned
\tag 4.3
$$
with $c'>0$, implying that (4.2) holds with $c=c_0-\delta $.

The given nonsmooth symbol $p$ can be approximated by smooth
symbols $p_k$ as described around (2.12). For  small
$\varepsilon _1,\varepsilon _2\in \,]0,c_0[\,$, we let $P_k=\Op(p_k)$ with symbol $p_k=\varrho
_k*p$, choosing $k$ so large that $|p_0(x,\xi )-p_{k,0}(x,\xi
)|\le \varepsilon _1|\xi |^{2a}$ for $|\xi |\ge 1$, and
$\|P-P_k\|_{\Cal L(H^a,H^{-a})}\le \varepsilon_1$, and we apply (4.3)
to $P_k$ with $c_0$ replaced by $c_0-\varepsilon _1$, $\delta
$ replaced by $\varepsilon _2$. This gives
$$
\aligned
\operatorname{Re}(Pu,u)&=\operatorname{Re}(P_ku,u)+\operatorname{Re}((P-P_k)u,u)
\\
&\ge ((c_0-\varepsilon _1)-\varepsilon _2)\|u\|^2_{H^a}-\beta '\|u\|^2_{L_2}-\|(P-P_k)u\|_{H^{-a}}\|u\|_{H^a}\\
&\ge (c_0-(2\varepsilon _1+\varepsilon _2))\|u\|^2_{H^a}-\beta '\|u\|^2_{L_2}. 
\endaligned
$$
Since any $\delta \in \,]0,c_0[\,$ can be written as $2\varepsilon _1+\varepsilon _2$,
this shows the assertion.\qed

\enddemo

Since $\dot H^{a}(\comega)\subset H^a(\rn)$, 
$r^+P$ defines in particular a continuous operator $\Cal P$ from $\dot
H^a(\comega) $ to $\ol H^{\,-a}(\Omega )$. Define the
Dirichlet realization $P_{D,2}$ as the operator  in
$L_2(\Omega )$ acting
like $r^+P$ with domain
$$
D(P_{D,2})= \{u\in \dot H^a(\comega)\mid r^+Pu\in
L_2(\Omega )\}. \tag4.4
$$

This can be viewed in a variational framework:  
Define the sesquilinear form
$$
s(u,v)=\int_{\Omega }Pu\,\bar v\,dx,\tag4.5
$$
first for $u,v\in C_0^\infty (\Omega )$, then extended by closure to
a continuous sesquilinear form on $\dot H^a(\comega)$.
Since $P$ is
strongly elliptic, the form is coercive in view of Lemma 4.1:
$$
\operatorname{Re}s(u,u)\ge c\|u\|^2_{\dot H^a(\comega)}-\beta \|u\|^2_{L_2(\Omega )}, \text{ with } c>0 \text{
and }\beta \in{\Bbb R} .\tag 4.6
$$
 Then the Lax-Milgram lemma (as recalled in e.g.\ \cite{G09, Sect.\
 12.4}) applies. For one thing $s(u,v)$ equals $\langle \Cal P
 u,v\rangle_{\ol H^{-a}(\Omega ), \dot H^a(\comega)}$, where $\Cal P\colon\dot
H^a(\comega)\to\ol H^{\,-a}(\Omega )$ acts like $r^+P$; moreover this induces
 the operator $P_{D,2}$ with domain  $\{u\in \dot
 H^a(\comega)\mid \Cal Pu\in
L_2(\Omega )\}$, the same as the domain described in (4.4). The
 inequality (4.6) holds for $u$ in the domain with $s(u,u)$ replaced by $(Pu,u)_{L_2(\Omega )}$.

When $P$ moreover has even symbol and $\tau >2a$,  we have from Theorem 3.2 that
$$
D(P_{D,2})=H^{a(2a)}(\comega),\tag4.7
$$
further described in (2.20)ff.\ and Remark 3.8.

When $\beta =0$ in (4.6),  $\Cal P\colon\dot
H^a(\comega)\to\ol H^{\,-a}(\Omega )$ is a homeomorphism with 
$$
\|\Cal Pu\|_{\ol H^{\,-a}(\Omega )}\ge c\|u\|_{\dot H^a(\comega)}, \tag4.8
$$
 and
$P_{D,2}$ has lower bound
$\inf\{\operatorname{Re}(P_{D,2}u,u)/\|u\|_{L_2}^2\mid u\in
D(P_{D,2})\setminus 0\}\ge c$ and is bijective from
$D(P_{D,2})$ to $L_2(\Omega )$. In general this
holds for $P+\beta I$ instead of $P$. Moreover, elementary estimates of the
numerical range $\nu (P_{D,2})=\{(P_{D,2}u,u)/\|u\|_{L_2}^2\mid u\in
D(P_{D,2})\setminus 0\}$ give (as in \cite{G09, Cor.\ 12.21})
that the spectrum and the numerical range are contained in a sectorial region
$$
M=\{\lambda \in \C\mid \operatorname{Re}\lambda \ge c-\beta ,
|\operatorname{Im}\lambda |\le Cc^{-1}(\operatorname{Re}\lambda +\beta ) \},\tag4.9
$$
where $C$ is a positive constant for which $|s(u,u)|\le
C\|u\|^2_{\dot H^a(\comega )}$ on $\dot H^a(\comega )$. In particular,
the resolvent $(P_{D,2}-\lambda )^{-1}$ exists for
$\lambda $ outside $M$, and its operator norm satisfies
$$
\|(P_{D,2}-\lambda )^{-1}\|_{\Cal L(L_2(\Omega ))}\le
C'\ang\lambda ^{-1}\text{ for }\operatorname{Re}\lambda \le
-\beta . \tag4.10             
$$

Since the
injection of $\dot H^a(\comega)$ into $L_2(\Omega )$ is compact, the spectrum of
$P_{D,2}$ is discrete, lying in $M$. We shall denote
$$
\Sigma = \text{ the spectrum of }P_{D,2}. \tag4.11
$$
 For $\lambda \in\Sigma $, $P_{D,2}-\lambda $ is a Fredholm operator, with index 0 since
the index depends continuously on $\lambda $.

Let us collect the outcome in a theorem:

\proclaim{Theorem 4.2} Assume Hypothesis {\rm 3.1}.

The Dirichlet realization $P_{D,2}$ in $L_2(\Omega )$
has domain {\rm (4.7)}. This equals $\dot H^{2a}(\comega)$ when
$a\in\,]0,\frac12[\,$; it is contained in
$\dot H^{1-\varepsilon }(\comega )$ when $a=\frac12$; and it is
contained in $\dot H^{2a}(\comega)+d^a\ol H^{a}(\Omega )$ and in $\dot
H^{a+\frac12-\varepsilon }(\comega)$ when $a\in
\,]\frac12,1[\,$ (locally if $\tau <1$).

With $c_0$ satisfying {\rm (4.1)}, there is for any $c\in \,]0,c_0[\,$ a number $\beta \in \R$ such that
$$
\operatorname{Re}(P_{D,2}u,u)_{L_2(\Omega )}\ge
c\|u\|^2_{\dot H^a(\comega)}-\beta \|u\|^2_{L_2(\Omega )}\text{ for }u\in D(P_{D,2}) .\tag 4.12
$$

The spectrum $\Sigma $ of $P_{D,2}$ is discrete and lies in $M$ {\rm (4.9)}, which also contains the numerical range.
If $\beta =0$, $P_{D,2}$ is a homeomorphism of
$H^{a(2a)}(\comega)$ onto $L_2(\Omega )$; more generally, $P_{D,2}-\lambda I$ has
this homeomorphism property when $\lambda \in \C \setminus \Sigma $,
and there is a resolvent estimate {\rm (4.10)}.   For $\lambda \in\Sigma $,  $P_{D,2}-\lambda I$ defines a Fredholm
operator with index zero from $H^{a(2a)}(\comega)$ to $L_2(\Omega )$.

\endproclaim

The resolvent set of 
$P_{D,2}$ is $\C\setminus \Sigma $. We denote by $N_\lambda $
the kernel of $P_{D,2}-\lambda I$; it is nontrivial only when $\lambda \in
\Sigma $. The $L_2$-adjoint $(P_{D,2})^*$ has as its spectrum
the conjugated set $\overline\Sigma $, and we denote the kernel of
$(P_{D,2})^*-\ol \lambda  I$ by $N'_{\ol \lambda }$. It is a cokernel of
$P_{D,2}-\lambda I$, in the sense that the range $(P_{D,2}-\lambda
I)D(P_{D,2})$ equals the orthogonal complement of  $N'_{\ol
\lambda }$ in $L_2(\Omega )$; this is the set of functions $f$ satisfying
$$
\int_\Omega f\ol\psi\, dx=0\text{ for }\psi \in N'_{\ol \lambda
}.\tag 4.13
$$
Here $\dim N'_{\ol\lambda }=\dim N_\lambda $.

The theorem can be used as in \cite{G18b} to establish solvability
properties of evolution problems
$Pu(x,t)+\partial_tu(x,t)=f(x,t)$; we shall follow this up
in Section 6. One can also ask about the asymptotic behavior of
eigenvalues; this is treated in a current work \cite{G22b} showing that the
expected asymptotic Weyl formula  holds for selfadjoint $P_{D,2}$, as
in \cite{G15b} for smooth cases.

Also for general $q\in \,]1,\infty [\,$, a Dirichlet realization
$P_{D,q}$ can be defined, namely the operator acting like $r^+P$ with
domain (cf.\ (3.3) and Remark 3.8) 
$$
D(P_{D,q})=\{u\in \dot H_q^a(\comega)\mid r^+Pu\in \ol L_q(\Omega )\}=H_q^{a(2a)}(\comega).\tag4.14
$$

\subsubhead 4.2 The regularity of eigenfunctions
\endsubsubhead

We shall now study the structure of the eigenfunctions  of the $L_q$-Dirichlet realizations of $P$; i.e.\ the nontrivial solutions of 
$$
r^+Pu=\lambda u ,
\;  u\in H_q^{a(2a)}( \comega),\tag4.15
$$
for $\lambda \in{\Bbb C}$; here we
use the regularity results for the homogeneous Dirichlet problem.

Smoothness properties of eigenfunctions  were found earlier in the $C^\infty $-setting in
\cite{G15b}, and we shall employ a similar strategy, as far as it goes
when the limited H\"olderness of the symbol and the domain are taken
into account.

In the analysis we
need an observation on the comparison of H\"older spaces, when powers
of the distance to the boundary $d_0(x)$ enter into the picture: When $\Omega $
is a $C^{1+\tau }$-domain with $\tau >0$,
then for $a,b>0$ with $a+b<1+\tau $, and $a,b,a+b\notin{\Bbb N}$, 
$$
\dot C^{a+b}(\comega)\subset d_0^a \dot C^{b
}(\comega ).\tag4.16
$$
This is
undoubtedly very well known (and enters in some form in many papers), but since we
have not been able to find an elementary reference, we include a
 proof in the Appendix, see Lemma A.5. The result is extended to  more
 general distance functions
 $d(x)$ in Lemma A.6.
For the $a$-transmission
 spaces, this implies:

\proclaim{Lemma 4.3} Let $\Omega $ be a $C^{1+\tau }$-domain, $\tau
>0$, and let $0<a<1$. There holds, for $a<t<\tau +a$, 
$$
C_*^{a(t)}(\comega)\subset \dot C^t(\comega)+d^a
e^+C^{t-a}(\comega)\subset d^a e^+C^{t-a}(\comega),\text{ when }
t,t-a\notin \N.\tag 4.17
$$
When $\tau \ge 1$, the general distance function $d$ can here be
replaced by the more precise function
$d_0(x)=\operatorname{dist}(x,\partial\Omega )$ (near $\partial\Omega
$). When $\tau <1$, the inclusions holds in a local sense, as in Remark
{\rm 2.1}.
\endproclaim

\demo{Proof} When $\tau \ge 1$, the first inclusion in (4.17) holds
globally, with distance function $d$ or $d_0$ at convenience; both are
$C^{1+\tau }$-functions. Then the second inclusion follows by
application of Lemmas A.5 resp.\ A.6 to $\dot C^t(\comega)$, showing that 
$$
\dot C^t(\comega)\subset d_0^a\dot C^{t-a}(\comega)\subset d_0^ae^+C^{t-a}(\comega),
$$
resp.\ the same formula with $d_0$ replaced by $d$.

When $\tau <1$, the inclusions hold in each of the local coordinate patches
used to describe $C_*^{a(t)}$ (as in Remark 2.1); here Lemma A.6 is
applied. \qed
\enddemo

In some of the formulations in the following, the extension by zero $e^+$
is tacitly understood.

Our result on the eigenfunctions is as follows:

\proclaim{Theorem 4.4} Assume Hypothesis {\rm 3.1}.

Let $P_{D,q}$ be the $L_q$ Dirichlet
realization, for some $q\in\,]1,\infty [\,$. 
The eigenfunctions of $P_{D,q}$ satisfy:

$1^\circ$ If $0$ is an eigenvalue of $P_{D,q}$, its associated
eigenfunctions $u_0$ are in
$H_p^{a(s+2a)}(\comega)$ for any $p\ge q $, any $s<\tau -2a$, 
hence also in
$ C_*^{a(\tau -\varepsilon )}(\comega)\subset d^aC^{\tau
-a-\varepsilon }(\comega)$ for  $\varepsilon >0$ (with $\tau
-a-\varepsilon, \tau -\varepsilon  \notin{\Bbb N}$).

$2^\circ$ For nonzero eigenvalues $\lambda $ there holds:
The
eigenfunctions $u_\lambda  $ of $P_{D,q} $ are in
$H_p^{a(t)}(\comega)$ for all $p\ge q$ and all $t\le 3a$ with $t<\tau
$. Hence they are in
$ C_*^{a(t )}(\comega)\subset d^aC^{t-a }(\comega)$ for
$t<\min\{3a,\tau \}$ (with $t-a,t  \notin{\Bbb N}$), and 
$$
u_\lambda \in d^aC^{\min\{2a, \tau -a\}-\varepsilon }(\comega).
\tag4.18
$$

The inclusions hold in a local sense (cf.\ Remark {\rm 2.1}) when
$\tau <1$. In all cases, $u_\lambda \in \dot C^a(\comega)$.

\endproclaim

\demo{Proof}  
When $\lambda  $ is an
eigenvalue, the associated eigenfunctions $u _\lambda $
are the nontrivial solutions of (4.15).

$1^\circ$. If $\lambda  =0$, $r^+Pu_0=0\in L_p(\Omega )$ for all $p\ge
 q$, so  Theorem 3.4 gives that $u_0\in \dot H_p^{a}(\comega)$
 for all $p\ge q$. Then furthermore, Theorem 3.2 gives that  $u_0\in H_p^{a(t)}(\comega)$
 for any $t<\tau
 $.

In view of (2.32), we then
also have that $u_0\in C_*^{a(\tau  -\varepsilon )}(\comega)$, any
$\varepsilon >0$. By Lemma 4.3,
$$
  u_0\in d^aC^{\tau -a-\varepsilon }(\comega),\tag4.19
$$
in a local form if $\tau <1$. 
This is applicable since $\tau -\varepsilon <\tau $, a fortiori $\tau -a-\varepsilon <\tau $.

$2^\circ$. Now consider a nonzero eigenvalue with eigenfunction $u_\lambda $.
Since $u_\lambda =\frac1\lambda r^+Pu$, we have from Theorem 3.2 that
$u_\lambda \in  H_{q_1}^{a(\min\{3a,\tau \}-\varepsilon )}(\comega)\subset H_{q_1}^{a(2a)}(\comega)$,
$q_1=q$. Using
the sequence constructed in Theorem 3.4, we find successively  for
$j=1,2,\dots$ that
$u_\lambda \in \dot H_{q_j}^a(\comega)$, hence $r^+Pu_\lambda =\lambda u_\lambda \in \dot
H_{q_j}^a(\comega)\subset \ol
H_{q_j}^a(\Omega)$, imply $u\in  H_{q_j}^{a(2a)}(\comega)\subset \dot
H_{q_{j+1}}^a(\comega)$. Thus $u_\lambda \in \dot
H_p^a(\comega)$ for all $p\ge q$.

Since $r^+Pu_\lambda =\lambda u_\lambda \in \dot H_p^a(\comega )\subset \ol H_p^a(\Omega
)$, we get more precisely from Theorem 3.2 that $u_\lambda \in
H_p^{a(t)}(\comega)$ for all $p\ge q$ and all $t\le 3a$ with $t<\tau $, and
consequently   by (2.32) and Lemma 4.3,
$$
u_\lambda \in C_*^{a(t)}(\comega)\subset
d^aC^{t-a}(\comega), \text{ for }t<\min\{3a,\tau \}, \;t, t-a\in \rp\setminus\N;
\tag 4.20$$
the inclusion holds in a local sense when $\tau <1$.
Hence $u_\lambda \in d^aC^{\min\{2a,\tau -a\}-\varepsilon }(\comega)$.
\qed

\enddemo

For smooth domains and symbols, it was shown in \cite{G15b} that $u_0\in
\E_a(\Omega )$ for $\lambda =0$, and $u_\lambda \in d^aC^{2a-\varepsilon
}(\comega)$ for $\lambda \ne 0$; this can be improved to
$d^aC^{2a}(\comega)$ when $a\ne \frac12$ by use of the precise formulas in
\cite{G14b} for how the operators act in
$C_*^s$-spaces.
As shown in \cite{G19},
the regularity
of $u_\lambda $ in the case $\lambda \ne 0$ cannot in general be lifted to
$d^aC^{2a+\delta }(\comega)$ with $\delta >0$ (by comparison of Taylor
expansions at $\partial\Omega $). We expect the same to be the
case in situations with finite smoothness.

\proclaim{Corollary 4.5} Assume Hypothesis {\rm 3.1}.

The Dirichlet realizations $P_{D,q}$,
$1<q<\infty $, have the same discrete set of eigenvalues and
eigenfunctions as $P_{D,2}$ for all $q\in \,]1,\infty [\,$.
\endproclaim

\demo{Proof} When $\lambda $ is an eigenvalue and $u_\lambda $ an
associated eigenfunction for some $q$, it also so for any $p<q$, since
$L_q(\Omega )\subset L_p(\Omega )$. For $p>q$, the regularity shown in Theorem
3.4 implies that $u_\lambda $ is an eigenfunction with the same
$\lambda $. We know that the set of eigenvalues in
case $q=2$ is discrete with finite-dimensional eigenspaces.
These eigenvalues and eigenfunctions
have the same role for all other $q$. \qed
\enddemo

These considerations show in particular that the operator in $L_q(\Omega )$,
$r^+P-\lambda \colon H_q^{a(2a)}(\comega)\to L_q(\Omega )$, has the same
finite-dimensional nullspace (kernel) $N_\lambda $ as in the case
$q=2$; it lies in $d^aC^{\min\{2a,\tau -a\}-\varepsilon
}(\comega)$.

For $q\ge 2$, this leads immediately to Fredholm properties of $P_{D,q}-\lambda $:

\proclaim{Corollary 4.6} Assume Hypothesis {\rm 3.1}, and let
$q\ge 2$. Consider
$$
P_{D,q}-\lambda \colon H_q^{a(2a)}(\comega)\to L_q(\Omega ).
$$
For $\lambda \in \C\setminus\Sigma $, this is a homeomorphism.
For $\lambda \in \Sigma $, this is a Fredholm operator with kernel
$N_\lambda $ and cokernel $N'_{\ol\lambda }$, in the sense that the
range consist of the function $f\in L_q(\Omega )$ satisfying {\rm
(4.13)}.
In particular, the spectrum of $P_{D,q}$ is $\Sigma $.
\endproclaim

\demo{Proof} Let $\lambda \in\Sigma $. Since $L_q(\Omega )\subset
L_2(\Omega )$, and $N_\lambda \subset \dot C^a(\comega)\subset L_q(\Omega )$ the nullspace of
$P_{D,q}-\lambda $ equals $N_\lambda $, as already noted. Moreover,
$N'_{\ol\lambda }\subset L_2(\Omega )\subset L_{q'}(\Omega )$. Since
$P_{D,2}-\lambda $ is surjective from $H^{a(2a)}(\comega)$ onto the functions
in  $L_2(\Omega )$ satisfying (4.13),  its restriction $P_{D,q
}-\lambda $ is surjective from $H_q^{a(2a)}(\comega)$ onto the functions
in  $L_q(\Omega )$ satisfying (4.13).

When $\lambda \in \C\setminus\Sigma $, we can argue in the same way,
replacing $N_\lambda $ and $N'_{\ol\lambda }$ by zero spaces.
\qed

\enddemo

The mapping properties can be lifted to $H_q^s$-spaces with higher
$s$, but when $\lambda \ne 0$, we must here take into account that the multiplication by
$\lambda $ is limited by the possibility to embed
$H_q^{a(s+2a)}(\comega)$ in $\ol H^s(\Omega )$.
\proclaim{Proposition 4.7} Assume Hypothesis {\rm 3.1} $1^\circ$,
 and let $\lambda \in{\Bbb C}$.

Let $-a\le s<\tau -2a$. Then $r^+P-\lambda $ maps continuously
$$
r^+P-\lambda \colon H_q^{a(s+2a)}(\comega)\to \ol H_q^{s'}(\Omega ),\text{
where }s'=\min\{s,a+\tfrac1q-\varepsilon \}, \text{ for small }\varepsilon >0.\tag4.21
$$

Assume moreover that Hypothesis {\rm 3.1} $2^\circ$ holds. 
Let $u\in \dot H_q^a(\comega)$
satisfy $(r^+P-\lambda )u\in 
\ol H_q^s(\Omega )$.
If $s<a+\frac1q$, then $u\in H_q^{a(s+2a)}(\comega)$. If $s\ge
a+\frac1q$,  then $u\in H_q^{a(3a+\frac1q-\varepsilon )}(\comega)$, any
$\varepsilon >0$. In particular,  $u\in H_q^{a(3a )}(\comega)$ if $s\ge a$.
\endproclaim

\demo{Proof} We have (3.1) for $r^+P$ alone. As for the multiplication
by $\lambda $, note that in view of (2.20) and (2.21),
$$
 H_q^{a(s+2a)}(\comega)\cases = \dot
H_q^{s+2a}(\comega)\text{ if }s+a<\frac1q\\
\subset\dot H_q^{a+\frac1q-\varepsilon }(\comega)\text{ if
}s+a\ge \tfrac1q.
\endcases\tag4.22
$$
Therefore the multiplication by $\lambda $ satisfies:
$$
\lambda \colon H_q^{a(s+2a)}(\comega)\to \cases  \dot
H_q^{s+2a}(\comega)\subset \ol H^s_q(\Omega )\text{ if }s+a<\frac1q\\
\dot H_q^{a+\frac1q-\varepsilon }(\comega)\text{ if
}s+a\ge \tfrac1q.
\endcases\tag4.23
$$
The combined operator $r^+P-\lambda $ maps into the space with the smallest exponent (4.21).

Now assume moreover that $P$ is strongly alliptic and $s\ge 0$. Since $u\in \dot
H_q^a(\comega)$,  $\lambda u\in \dot
H_q^a(\comega)$, so $r^+Pu\in \ol H_q^s(\Omega )+ \dot
H_q^a(\comega)$. Then we conclude from Theorem 3.2 that $u\in
H_q^{a(t+2a)}(\comega) $ with $t=\min\{s,a\}$. If $s\le a$, the proof
is complete. If $s>a$, we iterate the argument. Note first that for
$s>a+\frac1q-\varepsilon $, (4.23) at best gives $r^+Pu-\lambda u\in
\ol H^{a+\frac1q-\varepsilon }(\Omega )$ and hence $u\in
H^{a(a+\frac1q+2a-\varepsilon) }(\comega)=H^{a(3a+\frac1q-\varepsilon
)}(\comega)$; it cannot be lifted further in this way. 

Let $a<s<a+\frac1q-\varepsilon $, i.e., $s=a+\delta $ with
$0<\delta <\frac1q-\varepsilon $. An application of Theorem 3.2 gives 
$$u\in
H^{a(3a+\delta ) }(\comega)\subset\cases  \dot
H_q^{3a+\delta }(\comega)\text{ if }2a+\delta <\frac1q\\
\dot H_q^{a+\frac1q-\varepsilon }(\comega)\text{ if
}2a+\delta \ge \tfrac1q.
\endcases
$$
If $2a+\delta <\frac1q$, we iterate the argument to find 
$$u\in
H^{a(5a+\delta ) }(\comega)\subset\cases  \dot
H_q^{5a+\delta }(\comega)\text{ if }4a+\delta <\frac1q\\
\dot H_q^{a+\frac1q-\varepsilon }(\comega)\text{ if
}4a+\delta \ge \tfrac1q.
\endcases
$$
For very small values of $a$, further iterations may be needed, each
 providing a lift by $2a$. Eventually after $N$ steps, $2Na+\delta $ will surpass $\frac1q$, so
 that we can conclude $u\in H_q^{a(3a+\frac1q-\varepsilon )}(\comega)$.\qed

\enddemo

For a small $\varepsilon >0$, define 
$$
r_q=\min\{2a, a+\tfrac1q-\varepsilon \},\tag4.24
$$
then the outcome of Proposition 4.7 is that $r^+P-\lambda \colon H_q^{a(s+2a)}(\comega)\to \ol H_q^{s}(\Omega ) $ has the same
regularity properties as $r^+P$ when $s\in [0,r_q]$. Then Corollary
4.6 extends easily to $H_q^s$-spaces:

\proclaim{Theorem 4.8}  Assume Hypothesis {\rm 3.1},  let $s<\tau -2a$,
$q\ge 2$, $\lambda \in\C$, and if $\lambda \ne 0$ let $s\in [0,r_q]$. 
Consider
$$
r^+P-\lambda \colon H_q^{a(s+2a)}(\comega)\to \ol H_q^s(\Omega ).
$$
For $\lambda \in \C\setminus\Sigma $, this is a homeomorphism.
For $\lambda \in \Sigma $, this is a Fredholm operator with kernel
$N_\lambda $, and cokernel $N'_{\ol\lambda }$ in the sense that the
range consist of the function $f\in \ol H_q^s(\Omega )$ which satisfy {\rm
(4.13)}.
\endproclaim

\demo{Proof} This follows from Corollary 4.6 by restriction of
$P_{D,q}-\lambda $ to $ H_q^{a(s+2a)}(\comega)$. \qed
\enddemo

Corollary 4.6 and Theorem 4.8 will be extended to all  $q\in \,]1,\infty [\,$ in Theorems
4.16--4.18 below.

Note that the theorem shows existence and uniqueness, resp.\ Fredholm
solvability, for the Dirichlet problem
$$
Pu-\lambda u=f
\text{ in }\Omega ,\quad u=0\text{ in }\rn\setminus\Omega,
$$
with $f\in \ol H_q^s(\Omega )$ with $s\ge 0$, $q\ge 2$, $u$ given in
$\dot H^a(\comega)$. This has the following corollary for $f\in
C^s(\comega)$ when $\lambda =0$:

\proclaim{Corollary 4.9} Assume  Hypothesis {\rm 3.1}, and let $0\le
s<\tau -2a$. Consider the homogeneous Dirichlet problem {\rm (3.2)},
with $f\in C^s(\comega)$, and $u$ a priori assumed to be in $\dot H^a(\comega)$. 

If $0\notin \Sigma $, there is a unique solution $u$, which satisfies 
$$
u\in C_*^{a(s+2a-\varepsilon )}(\ol\Omega)\subset d^aC^{s+a-\varepsilon }(\comega),\tag4.25
$$
for small $\varepsilon >0$ (with noninteger $s+a-\varepsilon $).

If $0\in\Sigma $, there is a solution $u$, unique modulo $N_0$, when
$f$ satisfies {\rm (4.13)}; here $u$ satisfies {\rm (4.25)}.

\endproclaim

\demo{Proof} By (2.31), $f\in \ol H_q^{s-\varepsilon '}(\Omega )$ for
any $q\in \,]1,\infty [\,$, any $\varepsilon '>0$. We know from Corollary 4.6 that there is unique resp.\
Fredholm solvability, and $u$ lies in $H_q^{a(s-\varepsilon
'+2a)}(\comega)$. Letting $q\to\infty $, we find (4.25) (with a local
interpretation when $\tau <1$). \qed
\enddemo

\subsubhead 4.3 Realizations of the adjoint operator \endsubsubhead

When $q<2$, we need some additional information on the adjoint $P^*$,
to get useful regularity properties of $N'_{\ol\lambda }$.

The formal adjoint of $P$ on  $\rn$ (the pseudodifferential operator
$P^*$ in $\Cal S'(\rn)$ that satisfies $\ang {P^*u,\varphi }=\ang
{u,P\varphi }$ for all $\varphi \in \Cal S(\rn)$) is the operator $P^*=\Op(\overline
p(y,\xi ))$ in $y$-form (in the general concept of operators defined
from symbols $a(x,y,\xi )$ explained in \cite{AG21}). We here work
with a sesquilinear duality for consistency with the $L_2$-duality. $P^*$ is the sum of the $x$-form operator
$\overline P=\Op(\ol p(x,\xi ))$ and a difference operator $P'$
$$
P^*=\ol P+P',\tag4.26
$$
  where $\ol P$
has similar properties as $P$, and $P'$ is the remainder. Its mapping
properties were described in Marschall \cite{M88, Cor.\ 3.6} which has
the following consequence in our setting:

\proclaim{Proposition 4.10} Let $A=\Op(a(x,\xi ))$, where $a(x,\xi )\in C^\tau S^m(\rn\times\rn)$ for a
$\tau >0$ and an $m\ge 0$, let $\ol A=\Op (\ol a(x,\xi ))$ and let
$A'=A^*-\ol A$. Then $A'$ maps continuously, for $ \theta\in [0,1]$ with $\theta <\tau $,
$$
A'\colon H_q^{s+m-\theta }(\rn)\to H_q^s(\rn),\text{ when } -\tau +\theta <s<\tau .\tag 4.27
$$
\endproclaim

\demo{Proof} Corollary 3.6 in \cite{M88} states (with some relabeling
of parameters)
that when $a(x,\xi )$ is in $H^\sigma _QS^m_{1,0}(\rn\times\rn)$,
namely the
symbol space with estimates
$$
|a(x,\xi )|\le C\ang\xi ^m,\quad \|\partial_\xi ^\alpha a(\cdot,\xi
 )\|_{H_Q^\sigma (\rn)}\le C_\alpha \ang\xi ^{m-|\alpha |} ,\text{ all }\alpha ,
$$
 then for $0< p\le \infty $, $0\le\theta \le 1$, $n/Q+\theta <\sigma $, $0<Q\le \infty $
 and 
$$
n(\tfrac1Q +\tfrac1p-1)_+-\sigma +\theta <s<\sigma -n(\tfrac1Q-\tfrac1p)_+,
$$
the mapping $A'\colon H_p^{s+m-\theta }(\rn)\to H_p^s(\rn)$ is
bounded.
For a given $p$, we shall take $Q\in\rp$ so large that $(\tfrac1Q
+\tfrac1p-1)_+$ and $(\tfrac1Q-\tfrac1p)_+$ are 0. Since $C^\tau(\rn)
\subset H_Q^{\tau -\varepsilon }(\rn)$ (any small $\varepsilon >0$), the given symbol is in the
space $H^\sigma _QS^m_{1,0}(\rn\times\rn)$ with $\sigma =\tau
-\varepsilon $ and any $Q$. Then the asserted continuity
follows when $\sigma =\tau -\varepsilon $,  $0\le\theta < \sigma $ and
$\theta \le 1$, and $-\sigma  +\theta <s<\sigma
$.
This reduces to the mentioned conditions since $\varepsilon $ can be
taken arbitrarily small.  \qed
\enddemo

For our operator $P'$ we conclude

\proclaim{Corollary 4.11} Assume Hypothesis {\rm 3.1} $1^\circ$. The operator $P'=P^*-\ol P$
maps continuously
$$
\align
P'&\colon H_q^{s }(\rn)\to H_q^s(\rn),\text{ when }2a\le 1,\; -\tau
+2a <s<\tau ,\tag4.28\\
P'&\colon H_q^{s+2a-1 }(\rn) \to H_q^s(\rn),\text{ when }2a>1,\; -\tau
+1 <s<\tau ;\tag4.29
\endalign
$$
in particular, the mapping properties hold for $0\le s<\tau $.
\endproclaim

\demo{Proof} Let $m=2a$. If $2a\le 1$, we can take
$\theta =2a$ (which is $<\tau $) in (4.27); this shows (4.28). If $2a>1$, hence $\tau >1$,  (4.27) holds
with $\theta =1$, this shows (4.29).  \qed
\enddemo

A first consequence is that $P'$ maps $H_q^a(\rn)\to
L_q(\rn)$. This follows obviously from (4.28) when $2a\le 1$. When  $2a>1$, 
$a=s+2a-1$ when $s=1-a$, so $P'$ maps $H_q^a(\rn)$ to
$H_q^{1-a}(\rn)\subset L_q(\rn )$ by (4.29).  A fortiori, $P'$ maps
$H_q^a(\rn)$ to $H_q^{-a}(\rn)$.
Then 
since $\ol P$  has the same mapping properties as
$P$, $r^+P^*$ is
continuous from $\dot H_q^{a}(\comega)$ to $\ol H_q^{-a}(\Omega ) $,
and there holds, by extension by continuity from $u,v\in C_0^\infty (\Omega )$,
$$
\ang {r^+Pu,v}_{\ol H^{-a}_{q}(\Omega ), \dot H_{q'}^a(\comega)}=\ang
{u,r^+P^*v}_{\dot H_{q}^a(\comega),\ol H^{-a}_{q'}(\Omega )}, \text{ for
}u\in \dot H_q^a(\comega), v\in \dot H_{q'}^a(\comega),\tag4.30
$$
again with  sesquilinear dualities.

 The adjoint  of $P_{D,2}$ in $L_2(\Omega )$ by the Lax-Milgram Lemma
 (cf.\ e.g.\ \cite{G09, Sect.\
 12.4}) is the realization $(P^*)_{D,2}$ of $r^+P^*$ in $L_2(\Omega )$ with domain
$$
D((P^*)_{D,2})=\{u\in \dot H^a(\comega )\mid r^+P^*u\in L_2(\Omega ) \}.
$$
  Here  $ r^+\ol Pu\in L_2(\Omega )\iff u\in H^{a(2a)}(\comega)$; and
since $r^+P'$ maps $\dot H^a(\comega )$ to
$
L_2(\Omega )$ as seen above, hence maps also the subset $H^{a(2a)}(\comega)$
to $L_2(\Omega )$, it follows that
$D((P^*)_{D,2}) =H^{a(2a)}(\comega)$.

We can henceforth drop the
parentheses in the notation, noting that the adjoint of $P_{D,2}$ in $L_2(\Omega )$ is the
operator $P_{D,2}^*$ acting like $r^+P^*$ with domain
$$
D(P^*_{D,2})=H^{a(2a)}(\comega)=\{u\in \dot H^a(\comega )\mid r^+P^*u\in L_2(\Omega ) \}.\tag4.31
$$

\example{Remark 4.12}
If $P$ is $x$-independent,  $P^*$ equals $\ol P$, which behaves in exactly the
same way as $P$; typical examples of $x$-independent operators are $(-\Delta
)^a$ and $(m^2-\Delta )^a$. Note also that if an $x$-dependent  $P$ is known to be formally
selfadjoint, i.e.\  $P^*=P$, no further analysis of $P^*$ is needed (this holds
for instance for fractional powers of a formally selfadjoint strongly
elliptic differential operator). But 
in general, a nontrivial difference operator $P'$ will give some limitations in the
analysis of
regularity properties for $P^*$.
\endexample

We can now show the following 
restricted variant
of Theorem 3.2 for $P^*$:

\proclaim{Theorem 4.13} Assume Hypothesis {\rm 3.1} $1^\circ$, and let
$0\le s<\tau -2a$.

$1^\circ$ Let $2a\le 1$, and set $t_0=a+\frac1q$. Then $r^+P^*$ maps
continuously, for any small $\varepsilon >0$,
$$
r^+P^*\colon H_q^{a(s+2a)}(\comega)\to \cases \ol H_q^{s}(\Omega ),\text{
when }s<t_0,\\
 \ol H_q^{t_0-\varepsilon }(\Omega ),\text{
when }s\ge t_0.\endcases \tag4.32
$$
If moreover Hypothesis {\rm 3.1} $2^\circ$ holds, one has when $u\in \dot H_q^a(\comega)$:
$$
 r^+P^*u\in \ol H_q^{s}(\Omega ) \implies
\cases u\in H_q^{a(s+2a)}(\comega)\text{ when }s<t_0,\\
u\in H_q^{a(t_0-\varepsilon +2a)}(\comega)\text{ when }s\ge t_0.\endcases\tag4.33
$$
In particular, $r^+P^*u\in \ol H_q^a(\Omega )$ implies $u\in H_q^{a(3a)}(\comega)$.

$2^\circ$ Let $2a> 1$, and set $t_1=1-a+\frac1q$. Then $r^+P^*$ maps
continuously, for any small $\varepsilon >0$,
$$
r^+P^*\colon H_q^{a(s+2a)}(\comega)\to \cases \ol H_q^{s}(\Omega ),\text{
when }s<t_1,\\
 \ol H_q^{t_1-\varepsilon }(\Omega ),\text{
when }s\ge t_1.\endcases \tag4.34
$$
If moreover Hypothesis {\rm 3.1} $2^\circ$ holds, one has when $u\in \dot H_q^a(\comega)$:
$$
r^+P^*u\in \ol H_q^{s}(\Omega ) \implies
\cases u\in H_q^{a(s+2a)}(\comega)\text{ when }s<t_1,\\
u\in H_q^{a(t_1-\varepsilon +2a)}(\comega)\text{ when }s\ge t_1.\endcases\tag4.35$$
 In particular, $r^+P^*u\in \ol H_q^{1-a}(\Omega )$ implies $u\in H_q^{a(1+a)}(\comega)$.

\endproclaim

\demo{Proof}

$1^\circ$. {\it The case} $2a\le 1$. By (4.28), $P'$ preserves $H^s(\rn)$, so
$r^+P'\colon \dot H_q^s(\comega )\to \ol H_q^s(\Omega )$ for $0\le
s<\tau $,  acting similarly to the multiplication by $\lambda $
considered in Proposition 4.7. Since Theorem 3.2 applies to $\ol P$,
we can apply the proof of Proposition 4.7 to $\ol P+P'$, concluding
that under Hypothesis 3.1 $1^\circ$,
$$
r^+(\ol P+P') \colon H_q^{a(s+2a)}(\comega)\to  \cases \ol H_q^{s}(\Omega ),\text{
when }s<a+\frac1q,\\
 \ol H_q^{a+\frac1q-\varepsilon }(\Omega ),\text{
when }s\ge a+\frac1q,\endcases
$$
for small $\varepsilon >0$; this shows (4.32).

When moreover Hypothesis {\rm 3.1} $2^\circ$ holds, we find as in Proporsition 4.7
that for small $\varepsilon >0$,
$$
u\in \dot H_q^a(\comega), r^+P^*u\in \ol H_q^s(\Omega )\implies\cases
u\in H_q^{a(s+2a)}(\comega)\text{ when }s<a+\frac1q,\\
u\in H_q^{a(3a+\frac1q-\varepsilon )}(\comega)\text{ when }s\ge a+\frac1q;\endcases
$$
this shows (4.33).

$2^\circ$. {\it The case} $2a>1$. Here $r^+P'\colon \dot H_q^{s+2a-1}(\comega)\to
\ol H_q^s(\Omega )$ for $0\le s<\tau $. Now (cf.\ (2.20)--(2.21)) 
$$
H_q^{a(s+2a)}(\comega)\subset H_q^{a(s-1+2a)}(\comega)\cases =\dot
H_q^{s-1+2a}(\comega )\text{ if }s-1+a<\frac1q,\\
\subset \dot H_q^{a+\frac1q-\varepsilon }(\comega )\text{ if }s-1+a\ge
\frac1q.\endcases
$$
In the first case, the space is  in view of (4.29) mapped by $r^+P'$ into $\ol
H_q^s(\Omega )$; in the second case into   $\ol
H_q^{1-a+\frac1q-\varepsilon }(\Omega )$. This shows the
forward mapping property (4.34).

Now 
assume moreover that $P$ is strongly alliptic, so that the full Theorem 3.2
is valid for $\ol P$. Let $u\in \dot H_q^a(\comega)$ with $r^+P^*u\in
\ol H_q^s(\Omega )$.
Since $r^+P'\colon \dot H_q^a(\comega)\to \ol H_q^{1-a}(\Omega )$ by (4.29), we have
that
$$
r^+\ol Pu=r^+P^*u-r^+P'u\in \ol H_q^s(\Omega )+\ol H_q^{1-a}(\Omega ).\tag4.36
$$
If $s\le 1-a$, we conclude $u\in H_q^{a(s+2a)}(\comega)$ from the
regularity property of $\ol P$. If $s>1-a$, this
argument gives that $u\in H_q^{a(1+a)}(\comega)$. There is a small improvement of the
latter property:
We have that  $H_q^{a(1+a)}(\comega)\subset \dot
H_q^{a+\frac1q-\varepsilon }(\comega)$ by (2.21), hence $r^+P'u\in \ol
H_q^{\frac1q+1-a-\varepsilon }(\Omega )$. Insertion of this information in
(4.36) shows that if $s<\frac1q+1-a$, the regularity property of $r^+\ol P$
gives that $u\in H_q^{a(s+2a)}(\comega)$, and when $s\ge \frac1q+1-a$,
it gives that $u\in H_q^{a(\frac1q+1+a-\varepsilon
)}(\comega)$. This shows (4.35). \qed
\enddemo

For these results, recall that $s$ is also subject to the condition
$s<\tau -2a$; if $t_0$ or $t_1$ is larger, the range of possible $s$
may not include $t_0$ resp.\ $t_1$.

But $s=0$ is always included, and gives as a special case:

\proclaim{Corollary.4.14}  Assume Hypotheses {\rm 3.1} $1^\circ$.

$r^+P^*$ maps $H_q^{a(2a)}(\comega)$
continuously into $L_q(\Omega )$. When moreover  Hypotheses {\rm 3.1}
$2^\circ$
holds, we have: When $u\in \dot H_q^a(\comega)$
satisfies $r^+P^*u\in 
L_q(\Omega )$, then $u\in H_q^{a(2a)}(\comega)$.
\endproclaim

In the following, in indications of H\"older spaces with an $\varepsilon $,
it is understood that $\varepsilon $ is chosen so that integer
exponents are avoided.

For general $q$ we define the Dirichlet realization
$P^*_{D,q}$ of $P^*$ as the operator acting like $r^+P^*$ with domain
$$
D(P^*_{D,q})=\{u\in \dot H_q^a(\comega )\mid r^+P^*u\in L_q(\Omega ) \}=H^{a(2a)}_q(\comega);\tag4.37
$$
the last equality holds in view of Corollary 4.14. 

Again, we can show a regularity of eigenfunctions:

\proclaim{Theorem 4.15} Assume Hypothesis {\rm 3.1}. Let $u_\lambda
\in \dot H_q^a(\comega) $ be an eigenfunction of $r^+P^*$, i.e.\ a nontrivial solution of $r^+P^*u_\lambda =\lambda u_\lambda $.
Then 
$$
u_\lambda \in C_*^{a(\min\{2a, 1,\tau -a  \}+a-\varepsilon )}(\comega)\subset
d^a  C^{\min\{2a,1, \tau -a \}-\varepsilon }(\comega);
                                \tag4.38
$$
in a local sense if $\tau <1$.

As a consequence, the Dirichlet realizations $P^*_{D,q}$,
$1<q<\infty $, have the same discrete set of eigenvalues and
eigenfunctions as $P^*_{D,2}$ for all $q\in \,]1,\infty [\,$.
\endproclaim

\demo{Proof} 
Note first that in view of Corollary 4.14, the lifting of the
$q$-parameter in Theorem
3.4 can be performed also
with $P$ replaced by the present $P^*$. Therefore we can conclude that
$u_\lambda \in \dot H_p^a(\comega)$ for all $p<\infty $.
To use Theorem 4.13, we must combine the information given there with
the restriction $s<\tau -2a$.

In the case $2a\le 1$, we find from Theorem 4.13 $1^\circ$ that
$r^+P^*u_\lambda \in \dot H_p^a(\comega)$ implies
$$
u_\lambda \in H_p^{a(\min\{a, \tau -2a-\varepsilon , a+\frac1p-\varepsilon \}+2a)}(\comega).
$$
for small $\varepsilon >0$.
 In view  of (2.32) and Lemma 4.3, this implies
$$
u_\lambda \in C_*^{a(\min\{a, \tau -2a  \}+2a-\varepsilon )}(\comega)\subset
d^a  C^{\min\{2a, \tau -a \}-\varepsilon }(\comega),\tag4.39
$$
for small $\varepsilon >0$.

A better estimate for $\lambda =0$ cannot be expected (with the present
strategy), because of the
presence of $t_0-\varepsilon =a+\frac1p-\varepsilon $ in the estimate in
(4.33), no matter how large $s$ is. 

In the case $2a> 1$, we find from Theorem 4.13 $2^\circ$ that
$r^+P^*u_\lambda \in \dot H_p^a(\comega)$ implies
$$
u_\lambda \in H_p^{a(\min\{1-a, \tau -2a-\varepsilon , a+\frac1p-\varepsilon \}+2a)}(\comega).
$$
for small $\varepsilon >0$.
 In view  of (2.32) and Lemma 4.3, this implies
$$
u_\lambda \in C_*^{a(\min\{1-a, \tau -2a ,a \}+2a-\varepsilon
)}(\comega) = C_*^{a(\min\{1, \tau -a \}+a-\varepsilon
)}(\comega)
\subset
d^a  C^{\min\{1, \tau -a \}-\varepsilon }(\comega),\tag4.40
$$
for small $\varepsilon >0$. Now (4.39) and
(4.40) together imply (4.38).

The last statement follows as in Corollary 4.5.\qed

\enddemo

The theorem shows in particular that the cokernels $N'_{\ol\lambda }$ for $P_{D,2}-\lambda
$, $\lambda \in \Sigma $, have the regularity in (4.38).

\subsubhead 4.4 Fredholm operators in $L_q$ with a spectral parameter
\endsubsubhead

We now have the tools to deduce some spectral properties like those of $P_{D,2}$ for
the realizations $P_{D,q}$, besides what was shown for $q\ge 2$ in
Corollary 4.6 and Theorem 4.8.

It is well-known that $L_q(\Omega )$ is a reflexive Banach space, the
dual space identifying with $L_{q'}(\Omega )$, $\frac1{q'}=1-\frac1q$,
with a sesquilinear duality $\ang{f,g}_{L_q,L_{q'}}=\int_\Omega
f\,\bar g\, dx$. The adjoint $T^*$ in $ L_{q'}(\Omega )$ of a closed,
densely defined (unbounded) operator $T$ in $L_q(\Omega )$, has the
domain consisting of the functions $v$ for which there exist $v^*$
such that 
$$
\ang{Tu,v}_{L_{q},L_{q'} }=\ang {u,v^*}_{L_q,L_{q'}},\text{ all }u\in D(T);
$$
  then $v\in D(T^*)$ with $T^*v=v^*$ (uniquely). $T^*$ is likewise
  closed and densely defined, and $T^{**}=T$. This is usually deduced
  by a consideration of the graphs, where one also finds that $T$ is
  bijective if and only if $T^*$ is so. (Fredholm theory in reflexive Banach
  spaces is presented e.g.\ in Schechter \cite{S02}.)

For $q\ne 2$, the Dirichlet realizations of $r^+P$ and $r^+P^*$, as
unbounded operators in the reflexive Banach spaces $L_q(\Omega )$ resp.\
$L_{q'}(\Omega )$, will now be shown to be adjoints:

\proclaim{Theorem 4.16}  Assume Hypothesis {\rm 3.1}.

$1^\circ$ For each $1<q<\infty $, the adjoint of
$P_{D,q}$ in $L_q(\Omega )$ is the operator $P^*_{D,q'}$ in
$L_{q'}(\Omega )$.

$2^\circ$ For all $\lambda $ in the resolvent set $\C\setminus\Sigma $, $P_{D,q}-\lambda I$ is
bijective from $D(P_{D,q})$ to $L_q(\Omega )$. 
\endproclaim

\demo{Proof} In view of (4.30), there holds
$$
\ang{P_{D,q}u,v}_{L_{q},L_{q'} }=\ang
{u,P^*_{D,q'}v}_{L_q,L_{q'}},\text{ for all }u\in D(P_{D,q}), v\in D(P^*_{D,q'})
,
$$
so $D(P^*_{D,q'})\subset D((P_{D,q})^*)$, and $P^*_{D,q'}\subset
(P_{D,q})^*$ there.

Consider first the case where $0$ is not an eigenvalue of $P_{D,2}$, hence
not of $P^*_{D,2}$, nor of $P_{D,q}$ and $P^*_{D,q'}$ by Corollary
4.5 and Theorem 4.15.

Let $q<2$ so that  $q'>2$. We know that $P_{D,q}$ is injective, and
since $P_{D,q}$ extends $P_{D,2}$, its range contains $L_2(\Omega
)$. By Proposition 3.3, the range is closed, so since $L_2(\Omega )$
is dense in $L_q(\Omega )$, $P_{D,q}$ is surjective
from its domain to $L_q(\Omega )$; altogether it is bijective. Then the adjoint
$(P_{D,q})^*$ is likewise bijective, from its domain to $L_{q'}(\Omega )$. 

For $P^*_{D,q'}$, we know that it is injective (by Theorem 4.15), and
since $q'>2$ and $P^*_{D,2}$ is surjective, $P^*_{D,q'}$ is surjective
onto $L_{q'}(\Omega )$. Then the inclusion $P^*_{D,q'}\subset (P_{D,q})^*$
holds for two bijective operators, which implies that they are equal.

For $q>2$, hence $q'<2$, there is a similar proof where the roles of
$P$ and $P^*$ are reversed. This shows $1^\circ$.

In the case where $0$ is an eigenvalue, we can choose an arbitrary $\lambda _0\in
\C\setminus \Sigma $ and apply the above argumentation to the $L_q$
resp.\ $L_{q'}$
Dirichlet realizations of $P-\lambda _0$
and $(P-\lambda _0)^*=P^*-\ol \lambda _0$, showing that they are
adjoints and bijective. Here we use the regularity properties
of the shifted operator $P-\lambda _0I $  established in Proposition
4.7; there is a similar rule for $P^*-\ol\lambda _0I $.

Statement $2^\circ$ formulates the bijectiveness property of $P_{D,q}$
resp.\ $P_{D.q}-\lambda _0$ obtained in
the above proof. \qed

\enddemo

Moreover, the Fredholm properties of operators $P_{D,q}-\lambda $
shown in Corollary  4.6 for $q\ge 2$ can
now be extended to $q<2$:

\proclaim{Theorem 4.17}  Assume Hypothesis {\rm 3.1}.

For $\lambda \in \Sigma $, $P_{D,q}-\lambda $ is a Fredholm
operator from $H_q^{a(2a)}(\comega)$ to $L_q(\Omega )$ with index 
$0$. The kernel of $P_{D,q}-\lambda $ equals $N_\lambda $ (the kernel of
$P_{D,2}-\lambda $), and a
cokernel of $P_{D,q}-\lambda $ is
$N'_{\ol \lambda }$  (the cokernel of $P_{D,2}-\lambda $ equal to the
kernel of $P^*_{D,2}-\ol \lambda $), in the sense that the homogeneous
Dirichlet problem for $P-\lambda $,
$$
(P-\lambda )u=f \text{ in }\Omega ,\quad  u=0 \text{ in }\rn \setminus \comega, \tag4.41
$$
is solvable with $u\in H_q^{a(2a)}(\comega)$ when $f\in L_q(\Omega )$
satifies {\rm (4.13)} (i.e., $f$ is in the annihilator of
$N'_{\ol\lambda }$ in $L_q(\Omega )$).

Here $N_0\subset  d^aC^{\tau
-a-\varepsilon }(\comega)$, and $N_\lambda \subset   d^aC^{\min\{2a,
\tau -a\}-\varepsilon }(\comega)$, $N'_{\ol\lambda
}\subset d^a  C^{\min\{2a,1, \tau -a \}-\varepsilon }(\comega)$ for
general $\lambda \in\Sigma $ (locally when $\tau <1$); they are all
contained in $\dot C^a(\comega)$. For all $\lambda \in\Sigma $, $\operatorname{dim}N_\lambda =\operatorname{dim}N'_{\ol\lambda }$.
\endproclaim

\demo{Proof}
The last paragraph in the theorem repeats results from Theorems 4.2,  4.4
and 4.15. Since $\min\{2a,1,\tau -a\}>a$ and $d^a\in \dot C^a(\comega)$, all the eigenspaces are
contained in $\dot C^a(\comega)$.

For the main statement, consider for clarity first the case where $\lambda =0$, i.e., $0$ is
 an eigenvalue of $P_{D,2}$ with eigenspace $N_0$, and $P_{D,2}$ has a
cokernel $N_0'$ which is the eigenspace of $P^*_{D,2}$. Their
 regularity is recalled above.
By Corollary 4.5, $N_0$ is also the
eigenspace of $P_{D,q}$, and by Theorem 4.15, $N'_0$ is the
eigenspace of $P^*_{D,q'}$. Since $P_{D,q}$
 and $P^*_{D,q'}$ are adjoints, $N_0'$ is a cokernel of
 $P_{D,q}$. Thus $P_{D,q}$ is Fredholm with kernel $N_0$ and cokernel
 $N'_0$.

If $\lambda \ne 0$, we apply the same arguments to $(P-\lambda
)_{D,q}=P_{D,q}-\lambda $ and its adjoint $(P^*-\ol\lambda
)_{D,q'}=P^*_{D,q'}-\ol\lambda $, which in view of Proposition 4.7
and a similar rule for $P^*-\ol\lambda $ have sufficient regularity
properties to allow this.\qed

\enddemo

The resolvent problem (4.41) was discussed for selfadjoint operators like  $P=(-\Delta )^a$ by Chan, Gomez-Castro and Vazquez 
\cite{CGV21} with $f$ given (primarily) in a weighted  $L_1$-space; and a
Fredholm property was deduced from the knowledge of the solution
operator for $\lambda =0$, as an integral operator with kernel
$\Cal G(x,y)$ (Green's kernel) satisfying explicit estimates. The above
results develop such knowledge for a large class of  operators
(containing $(-\Delta )^a$) which can be $x$-dependent
and nonselfadjoint, with precise mapping properties in $L_q$ Sobolev
spaces for $1<q<\infty $.

We can also draw consequences for data in spaces with higher
regularity:

\proclaim{Theorem 4.18}  Assume Hypothesis {\rm 3.1}.

Consider the problem {\rm (4.41)} with $f$
given in $\ol H_q^s(\Omega )$ for some $0\le s<\tau -2a$. For a small
$\varepsilon >0$, let $r_q=\min\{2a, a+\frac1q-\varepsilon \}$.

If $\lambda
\notin\Sigma $, and $s\le r_q$ if $\lambda \ne 0$, {\rm (4.41)} is uniquely solvable with solution in
$H_q^{a(s+2a)}(\comega )$, and the solution operator, a restriction of $(P_{D,q}-\lambda
)^{-1}$, is continuous from $\ol H_q^s(\Omega )$ to
$H_q^{a(s+2a)}(\comega)$.

If $\lambda \in \Sigma $, and $s\le \min\{a,r_q\}$ if $\lambda \ne
0$, then  the problem {\rm (4.41)} is Fredholm
solvable from  $\ol H_q^s(\Omega )$ to
$H_q^{a(s+2a)}(\comega)$, in the sense that a solution $u\in
H_q^{a(s+2a)}(\comega) $ exists for
$f\in \ol H_q^s(\Omega )$
satisfying {\rm (4.13)}, and is unique modulo $N_\lambda $.
\endproclaim

\demo{Proof} For $\lambda =0$, the statements follow from Theorem 3.2
together with the added knowledge on bijectiveness or Fredholm solvability shown
in Theorems 4.16 and 4.17.

Now let $\lambda \ne 0$. We have from (4.22)
$$
 H_q^{a(2a)}(\comega)\cases = \dot
H_q^{2a}(\comega)\text{ if }a<\frac1q\\
\subset\dot H_q^{a+\frac1q-\varepsilon }(\comega)\text{ if
}a\ge \tfrac1q,
\endcases\tag4.42
$$
so that $H_q^{a(2a)}(\comega)\subset \dot H_q^s(\Omega )\subset \ol
H_q^s(\Omega )$ for $s\in [0,r_q]$. When $u$ solves (4.41) and $f\in \ol
H_q^s(\Omega )$ for some $s\in [0,r_q]$, then since $u\in
H_q^{a(2a)}(\comega) $ by Theorem 4.17, 
$$
r^+Pu=f+\lambda u\in \ol H_q^s(\Omega ),
$$
and it follows from Theorem 3.2 that $u\in
H_q^{a(s+2a)}(\comega)$.

When $\lambda \in \C\setminus \Sigma $, the
mapping from $f$ to $u$ is bijective as a consequence of Theorem 4.16,
hence continuous (by the closed graph theorem) since the inverse is so
by Theorem 3.2.

When $\lambda \in \Sigma $, we need the condition $s\le a$ to have
$N_\lambda \subset H_q^{a(s+2a)}(\comega)$, so that the Fredholm
solvability can be inferred from the property known from Theorem 4.17
for $s=0$. \qed

\enddemo

\example{Remark 4.19} In continuation of Remark 4.12, let us
underline that in cases where $P^*$ is as smooth as $P$, the
results hold with as good estimates for $P^*$ as for $P$. This goes
for $x$-independent operators (as for example $(m^2-\Delta )^a$,
$m\in{\Bbb R}$) and for 
selfadjoint cases as mentioned in Remark 4.12.

\endexample

\example{Example 4.20} In the case where $P$ has an $x$-independent
symbol $p(\xi )$ that is homogeneous of degree $2a$, even and positive
for $\xi \ne 0$, there is for any bounded set $\Omega \subset \rn$ a Poincar\'e inequality
$(Pu,u)\ge c\|u\|^2_{L_2(\Omega )}$ for $u\in \dot H^a(\comega)$, cf.\ e.g.\ the survey
of Ros-Oton \cite{R16, (3.4)ff.}. Then $P_{D,2}$ is
bijective. Moreover, the $L_q$-realizations $P_{D,q}$ are bijective
from $D(P_{D,q})=H_q^{a(2a)}(\comega)$ to $L_q(\Omega )$ when
$\Omega $ is $C^{1+\tau }$, for all $1<q<\infty $, by Corollary 4.5.
\endexample

\subhead 5. Solvability of the nonhomogeneous Dirichlet problem, spectral
analysis of  ``large'' solutions
\endsubhead

For the nonhomogenous Dirichlet problem, we have to go out in the
larger  spaces $H_q^{(a-1)(s+2a)}(\comega)$; here we get
solvability results by use of the  results in the homogeneous case for
$\lambda =0$, combined with Theorem 2.3.

\proclaim{Theorem 5.1}  Assume Hypothesis {\rm 3.1}  with
$\tau >2a+1$. Let  $0\le s<\tau -2a-1$.

Consider the nonhomogeneous Dirichlet problem {\rm (3.15)}, recalled here:
 $$
 \aligned
     Pu& = f \text { in }\Omega ,\\
     u&=0 \text { in }\rn\setminus \Omega ,\\
     \gamma _0^{a-1}u&=\varphi \text{ on }\partial\Omega ;
     \endaligned 
$$
with $f$
given in $\ol H^s_q(\Omega )$, $\varphi $ given in $
B_q^{s+a+1/q'}(\partial\Omega )$, and the solution being sought in 
$H_q^{(a-1)(s+2a)}(\comega)$.

If $0
\notin\Sigma $, it is uniquely solvable, and the solution operator,
given by the formula {\rm (5.3)}  below, is continuous from $\ol H^s_q(\Omega
)\times B_q^{s+a+1/q'}(\partial\Omega )$ to
$H_q^{(a-1)(s+2a)}(\comega)$.

If $0\in \Sigma $, it is Fredholm solvable, in the sense that a
solution $u$ exists, unique modulo $N_0$, when
$f'=f-r^+PK^{a-1}_{(0)}\varphi $ plays the role of $f$ in {\rm (4.13)}
for $\psi \in N_0'$.
\endproclaim

\demo{Proof}  Note that since $\tau >2a+1$, the case $\tau <1$ does
not occur here.

According to Theorem 2.3, there exists a right inverse $K^{a-1}_{(0)}$  of $\gamma _0^{a-1}$, mapping 
$$
K^{a-1}_{(0)}\colon  B_q^{s+a+1/q'}(\partial \Omega ) \to H_q^{(a-1) (s+2a)}(\comega ),
$$
for $-a-1/q'<s<\tau -2a-1$. 

Set $v=K^{a-1}_{(0)}\varphi $, then $u\in H_q^{(a-1)(s+2a)}(\comega)$ solves the given problem  if
and only if $w=u-v$ solves 
 $$
 \aligned
    Pw& = f -Pv \text { in }\Omega ,\\
     w&=0 \text { in }\rn\setminus \Omega ,\\
     \gamma _0^{a-1}w&=0 \text{ on }\partial\Omega .
     \endaligned \tag5.1
$$
Here $r^+Pv\in \ol H^s(\Omega )$, and by the last statement in Theorem
2.3,
$$
w\in \{u\in H_q^{(a-1)(s+2a)}(\comega)|\mid\gamma
_0^{a-1}u=0\}=H_q^{a(s+2a)}(\comega). \tag5.2
$$
so (5.1) is a homogeneous Dirichlet problem with right-hand side  $f-r^+Pv\in \ol
H^s(\Omega )$; $w$ being sought in $H_q^{a(s+2a)}(\comega)$. Here
Theorems 4.16 and 4.17 apply:

If $0\notin\Sigma $, the problem is uniquely solvable, with
$$
u=P_{D,q}^{-1}(f-r^+Pv)=P_{D,q}^{-1}f-P_{D,q}^{-1}r^+PK^{a-1}_{(0)}\varphi .\tag5.3
$$
If $0\in\Sigma
$, the problem is Fredholm solvable, in the sense that a solution $w$
exists if and only if $f'=f-Pv$ plays the role of $f$ in (4.13) with
$\psi \in N_0'$, and it is unique modulo $N_0$. This implies the
statement in the theorem for $u=w+v$. \qed

\enddemo

Observe that this is an existence-and-uniqueness theorem (resp.\
 Fredholm theorem), which completes the regularity result
Theorem 3.7 shown in Section 3. We can also complete Corollary 3.9
 with solvability in
 H\"older spaces:

\proclaim{Corollary 5.2} Assume Hypothesis  {\rm 3.1}  with
$\tau >2a+1$, and let  $0\le s<\tau -2a-1$. Let $f\in C^s(\comega)$ and $\varphi
\in C^{s+a+1}(\partial\Omega  ) $.

With $u$ a priori assumed to be in $H_q^{(a-1)(a)}(\comega)$ for some
$q$, problem {\rm (3.15)} is uniquely solvable if $0\notin\Sigma $, and
uniquely solvable modulo $N_0$ when $f'=f-r^+PK_{(0)}^{a-1}\varphi $
satisfies {\rm (4.13)}; and the solution satisfies $
u\in
C_*^{(a-1)(s+2a-\varepsilon )}(\ol\Omega)\subset \dot
C^{s+2a-\varepsilon }(\comega)+d^{a-1}C^{s+a+1-\varepsilon }(\comega)
$ as in Corollary {\rm 3.9}. 

\endproclaim

\demo{Proof} The statement follows by embedding the H\"older spaces in
$H^t_q$-spaces and applying Theorem 5.1, letting $q\to\infty $ as in the proof of Corollary 3.9.\qed
\enddemo

Next, we shall consider  nonhomogeneous problems with a spectral parameter $\lambda $ subtracted from
$P$. A study of such problems was initiated by Chan, Gomez-Castro and Vazquez  in \cite{CGV21}.
Since
the solutions generally blow up at the boundary (when $u/d^{a-1}$ has
a nonzero boundary value),
it is more demanding than in the homogeneous case to have $r^+Pu$ and $\lambda
u$ lying  in the same space. \cite{CGV21} handles this (for operators like
$P=(-\Delta )^a$) by considering $P$ in a weighted $L_1$-space. In our
treatment, we have $P$ defined on $H_q^{(a-1)(s+2a)}(\comega)$ and can get
results when this is contained in $L_q(\Omega )$.

\proclaim{Lemma 5.3} When $q<(1-a)^{-1}$, then 
$$
H_q^{(a-1)(s)}(\comega)\subset L_q(\Omega )\text{ for }s\ge 0.\tag5.4
$$

\endproclaim

\demo{Proof} Denote $a-1+1/q=t$; the hypothesis means that
$t>0$. Since the spaces $H_q^{(a-1)(s)}(\comega)$ dcrease with
increasing $s$, it suffices to show the statement for $0\le
s<t$. Working in local coordinates, we have for such $s$ (cf.\ (2.17))
$$
\aligned
H_q^{(a-1)(s)}(\crnp)&=\Xi _+^{-a+1}e^+\ol H_q^{s-a+1}(\rnp)
=\Xi _+^{-a+1}e^+\ol H_q^{s-t+1/q}(\rnp)=\Xi _+^{-a+1}\dot H_q^{s-t+1/q}(\crnp)\\
&= \dot
H_q^{s-t+1/q-a+1 }(\crnp) = \dot H_q^{s}(\crnp)\subset
L_q(\rnp).\qquad \square
\endaligned
$$
\enddemo

Thus for any given  $a\in \,]0,1[\,$, the inclusion (5.4) holds when $q$ is sufficiently low.
Note that for $q=2$,
$$
H^{(a-1)(s)}(\comega)\subset L_2(\Omega ),\text{ when }a>\tfrac12,
s\ge 0.\tag5.5
$$

\proclaim{Theorem 5.4} Assume Hypothesis {\rm 3.1}  with
$\tau >2a+1$. Assume moreover $q<(1-a)^{-1}$.
Consider the problem
 $$
 \aligned
     Pu-\lambda u& = f \text { in }\Omega ,\\
     u&=0 \text { in }\rn\setminus \Omega ,\\
     \gamma _0^{a-1}u&=\varphi \text{ on }\partial\Omega ,
     \endaligned \tag5.6
$$
with $f$
given in $L_q(\Omega )$, $\varphi $ given in $
B_q^{a+1/q'}(\partial\Omega )$, and the solution being sought in 
$H_q^{(a-1)(2a)}(\comega)$.

If $\lambda
\notin\Sigma $, it is uniquely solvable, and the solution operator,
given by the formula {\rm (5.8)}  below, is continuous from $L_q(\Omega
)\times B_q^{a+1/q'}(\partial\Omega )$ to
$H_q^{(a-1)(2a)}(\comega)$.

If $\lambda \in\Sigma $, it is Fredholm solvable,  in the sense that
there is a 
solution $u\in H_q^{(a-1)(2a)}(\comega)$, unique modulo $N_\lambda $, when
$f'=f-(r^+P-\lambda )K^{a-1}_{(0)}\varphi $ satisfies {\rm
(4.13)}. The solution operator is continuous from the closed subset of $L_q(\Omega
)\times B_q^{a+1/q'}(\partial\Omega )$ of pairs $\{f,\varphi \}$ such
that $f'=f-(r^+P-\lambda )K^{a-1}_{(0)}\varphi $ satisfies {\rm
(4.13)}, to
$H_q^{(a-1)(2a)}(\comega)$.

\endproclaim

\demo{Proof} Using the right inverse $K^{a-1}_{(0)}$ of $\gamma
_0^{a-1}$ recalled in the preceding proof,
we set $v=K^{a-1}_{(0)}\varphi $; then $u$ solves the problem (5.6) if
and only if $w=u-v$ solves 
 $$
 \aligned
    Pw-\lambda w& = f -(P-\lambda )v \text { in }\Omega ,\\
     w&=0 \text { in }\rn\setminus \Omega ,\\
     \gamma _0^{a-1}w&=0 \text{ on }\partial\Omega .
     \endaligned \tag5.7
$$
Here $\lambda v\in L_q(\Omega )$ by Lemma 5.3, so $ f -(r^+P-\lambda
)v\in L_q(\Omega )$. Problem (5.7) is in fact a homogeneous Dirichlet problem,
so by Theorem 4.16 it has the unique solution  $w=(P_{D,q}-\lambda )^{-1}(f-(r^+P-\lambda )v )$,
when $\lambda \notin \Sigma $. The solution of (5.6)
is then
$$
\aligned
u&=(P_{D,q}-\lambda )^{-1}(f-(r^+P-\lambda )v )+v \\
&=(P_{D,q}-\lambda
)^{-1}f+(1-(P_{D,q}-\lambda )^{-1}(r^+P-\lambda ))K^{a-1}_{(0)}\varphi.
\endaligned \tag5.8
$$

When $\lambda \in\Sigma $, we apply the Fredholm solvability of the
homogeneous Dirichlet problem (5.7) shown in Theorem 4.17.
\qed\enddemo

The theorem can be extended to slightly more smooth data with $s>0$,
subject to the condition $a-\frac1{q'}>s$.

In the discussion of the problem (5.6)  in \cite{CGV21} (for
operators like $P=(-\Delta )^a$),
it is called an {\it eigenvalue problem}. For $f\ne0$, we see it more
as a {\it resolvent problem} (where the solution operator may be useful e.g.\ in
associated evolution problems). 

However for $f=0$, (5.6) can certainly be regarded as an eigenvalue problem,
where solutions are sought that satisfy a fixed nonhomogeneous boundary
condition. Here \cite{CGV21} proved the existence of a nontrivial
solution when $\varphi $ is a continuous function $\ne 0$, by use of a resolvent for the
homogeneous Dirichlet problem (corresponding
to our $(P_{D,q}-\lambda )^{-1}$), acting in weighted $L_1$-spaces.

It is a major point in \cite{CGV21} that these
``eigenfunctions'' blow up at the boundary (like $d^{a-1}$),
being  ``large solutions''.

We agree of course that it is striking, that these solutions are
generally unbounded at $\partial\Omega $, but we also find that it is natural,
in view of the systematic point of view on how to define a nonhomogeneous local
Dirichlet condition, as presented at the start of this paper.

\example{Example 5.5} Here are some more details on the occurring
spaces in the case $q=2$, $a>\frac12$,  cf.\
(5.5). The boundary space is here $H^{a+\frac12}(\partial\Omega )$ and
the solution space is $H^{(a-1)(2a)}(\comega)$, which is a certain subspace of
$L_2(\Omega )$.
More technically, it is described in Remark 3.8 by:
$$
\aligned
H^{(a-1)(2a)}(\comega)&=H^{a(2a)}(\comega)\,\dot+ \,K^{a-1}_{(0)}H^{a+\frac12}(\partial\Omega
),\text{ where}\\
H^{a(2a)}(\comega)&=\dot H^{2a}(\comega)\, \dot+ \,K^{a}_{(0)}H^{a-\frac12}(\partial\Omega
);
\endaligned\tag5.9
$$
the operators $K^{a-1}_{(0)}$ and $K^{a}_{(0)}$ provide factors
$d^{a-1}$ resp.\ $d^a$.

\endexample

\subhead 6. Resolvent estimates and evolution problems \endsubhead

Solutions of evolution problems (parabolic problems) for $P$  with homogeneous boundary conditions were constructed in \cite{G18a,G18b} in
smooth settings, and we can now extend those results to the present
cases of operators $P$ and
domains $\Omega $ with limited smoothness. Moreover, we can introduce
completely new
results on evolution problems with
nonhomogeneous boundary conditions.

\subsubhead 6.1 Results for $q=2$ \endsubsubhead

Denote, for $\lambda \in \C\setminus \Sigma $, 
$$
(P_{D,2}-\lambda )^{-1}=R_\lambda ; \tag6.1
$$
it is the solution operator for the homogeneous Dirichlet problem for $r^+P-\lambda $
 in $L_2(\Omega )$. By Theorem 4.2, it is defined in particular for
 $\lambda $ in the
 complement of $M$
 (4.9). 

As a special case of (4.22),
$$
H^{a(2a)}(\comega)\cases =\dot H^{2a}(\comega)\text{ when
}0<a<\frac12,\\
\subset \dot H^{a+\frac12-\varepsilon }(\comega)\text{ when
}\frac12\le a<1,\endcases\tag6.2
$$
 any $\varepsilon >0$. Define
$$
r=\min\{2a,a+\tfrac12-\varepsilon \},\tag6.3
$$
for a small $\varepsilon \in \,]0, a+\tfrac12[\,$; then
$D(P_{D,2})\subset \dot H^r(\comega)\subset \ol
H^r(\Omega )$. (Here $r=r_2$ in (4.24).)

We know from Theorem 4.8 that $R_\lambda $ by resctriction defines a homeomorphism from $\ol
H^s(\Omega
)$ to
$H^{a(s+2a)}(\comega)$ for $\lambda $ in the
resolvent set $\C\setminus \Sigma $ including
$\C\setminus M$, when $s\in [0,r]$, $0\le s<\tau -2a$.

For the treatment of evolution problems we need norm estimates of $R_\lambda $ that are uniform in  $\lambda
$. To start with, there is the estimate (4.10); by the proof of \cite{G18b,
Th.\ 5.8} it can be supplied with estimates in spaces of higher
regularity.
Recall the notation from \cite{G18b} for a general operator $A$ in $L_2(\Omega )$:
$$
D_s(A)=\{u\in D(A)\mid Au\in \ol H^s(\Omega )\}\text{ for }s\ge 0;\tag6.4
$$
it will be applied to $A=P_{D,2}$. The space equals $H^{a(2a+s)}(\comega)$ when $s$ is as in  Theorem 3.2. 

\proclaim{Theorem 6.1} Assume Hypothesis {\rm 3.1}, and let
$s\in [0,r]$, $s<\tau -2a$.
 Then for $\lambda \notin M$, the resolvent $R_\lambda =(P_{D,2}-\lambda )^{-1}$ maps continuously
$$
R_\lambda :\ol H^s(\Omega )\simto H^{a(s+2a)}(\comega),\tag6.5
$$
satisfying the estimates when $\operatorname{Re}\lambda \le -\beta $:
$$
\aligned
\|R_\lambda f\|_{D_0(P_{D,2})} +\ang\lambda \|R_\lambda
f\|_{L_2(\Omega )} &\le C_0\| f\|_{L_2(\Omega )} ,\\
\|R_\lambda f\|_{D_s(P_{D,2})} +\ang\lambda ^{j+1}\|R_\lambda
f\|_{L_2(\Omega )} &\le C_j(\| f\|_{\ol H^s(\Omega )} +\ang\lambda
^{j}\| f\|_{L_2(\Omega  )})\text{ for }j\in\N.
\endaligned\tag6.6
$$
\endproclaim

\demo{Proof} The mapping property (6.5) is known from Theorem 4.8 (as
noted above).

The first (well-known) estimate in (6.6) is obtained by writing (4.10) as 
$$
\ang\lambda \|R_\lambda f\|_{L_2(\Omega )}\le C'\|f\|_{L_2(\Omega )},
$$
and supplying it with the observation
using that $P_{D,2}=(P_{D,2}-\lambda )+\lambda $,
$$
\|R_\lambda f\|_{D_0(P_{D,2})}\le c_1\|P_{D,2}R_\lambda f\|_{L_2}\le
c_1(\|f\|_{L_2}+\|\lambda R_\lambda f\|_{L_2})\le c_1(\|f\|_{L_2}+C'\| f\|_{L_2}).
$$

The second estimate in (6.6) is shown in \cite{G18b}, proof of Th.\ 6.8 (see in
particular formula (5.33) there), where it is seen how the result follows
from a combination of the regularity results
and numerical range estimates for the homogeneous Dirichlet
problem. Since these prerequisites hold for $P_{D,2}$ under the present
hypotheses on $\tau $ and $s$, the conclusion follows. The constant
denoted $\xi _0$ there equals $\beta $ in Theorem 4.2 here. \qed
\enddemo

One can possibly extend the second estimate in (6.6) to allow replacement of
$j\in \N$ by $j \in \rp$, but we think a nonzero $j $ is
needed if one wants to have estimates with $s>0$.

The estimates were used in \cite{G18b} to get  regularity estimates for
the solutions of the evolution problem (for some $T\in\rp$)
$$\aligned
Pu+\partial_tu&=f\text{ on }\Omega \times I ,\quad I=\,]0,T[\,,\\
u&=0\text{ on }(\Bbb R^n\setminus\Omega )\times I, \\
u|_{t=0}&=0;
\endaligned \tag6.7
$$
with zero initial value and homogeneous boundary condition.
Namely, as accounted for in the proofs of \cite{G18b,Th.\ 5.6 and
5.8}, there holds:

\proclaim{Theorem 6.2} Assume Hypothesis {\rm 3.1}, and let
$s\in [0,r]$, $s<\tau -2a$. Solutions
of {\rm (6.7)} are searched for $u$ in $ L_2(I;H^{a(2a)}(\comega ))$.

$1^\circ$ For $f$ given in $L_2(\Omega \times I)$, there is a
  unique solution $u$ of {\rm (6.7)} satisfying
  $$
u\in L_2(I;H^{a(2a)}(\comega ))\cap \ol H^1(I; L_2(\Omega ));\tag6.8
$$
moreover, $u\in \ol C^0(I;L_2(\Omega ))$.

$2^\circ$ Let $0\le s\le r$. If  $f\in L_2(I;\ol H^{s}(\Omega ))\cap \ol H^1(I;L_2(\Omega ))$, with $f|_{t=0}=0$, then the
solution of {\rm (6.7)} satisfies
$$
u\in L_2(I; H^{a(2a+s)}(\comega ))\cap \ol H^2(I;L_2(\Omega )).\tag6.9
$$

$3^\circ$ For any integer $j\ge 2$, if
$
f\in L_2(I;\ol H^{s}(\Omega
))\cap \ol H^j(I;L_2(\Omega ))$ with $\partial_t^lf|_{t=0}=0$ for
$l<j$, then
$$
u \in L_2(I; H^{a(2a+s)}(\comega ))\cap \ol H^{j+1}(I;L_2(\Omega )).\tag6.10
$$
It follows in particular that
$$
 f\in \bigcap_l \ol H^j(I;\ol H^s(\Omega
 )),\;\partial_t^lf|_{t=0}=0\text{ for } l\in{\Bbb N}_0\implies
u \in \bigcap_l\ol H^{j}(I; H^{a(2a+s)}(\comega )).\tag6.11
$$
\endproclaim

\demo{Proof} Statement $1^\circ$ follows as in \cite{G18b,Th.\
5.6}. Statements $2^\circ$ and $3^\circ$ are shown in \cite{G18b,Th.\
5.8} to follow from the second estimate in (6.6) by use of the
abstract result of Lions and Magenes \cite{LM68}
quoted as \cite{G18b, Th.\ 5.7}.\qed
\enddemo

Observe that the results allow high regularity in $t$, but that the
regularity in $x$ is at most achieved up to $s<\frac32$ (since
$r<\frac32$ for all $a\in \,]0,1[\,$). We think that this is not just due to the
method; there will in general be upper bounds on the $x$-regularity
for general data, as studied in detail in H\"older spaces in
\cite{G19}.

We can now also derive some results where nonhomogeneous boundary
contitions are included. Assume $\tau >1+2a$. Recall from Theorem 2.3
that $\gamma ^{a-1}_0$ has a continuous right inverse $K^{a-1}_{(0)}$,
 mapping 
$$
\aligned
&K^{a-1}_{(0)}\colon  H^{s-a+\frac12}(\partial \Omega ) \to H^{(a-1)
(s)}(\comega ), \text{ for }a-\tfrac12<s<\tau +a-1;\\
&\text{in particular, }K^{a-1}_{(0)}\colon  H^{a+\frac12}(\partial
\Omega )
\to H^{(a-1)(2a)}(\comega ).
\endaligned\tag6.12
$$
Now we need to assume $a>\frac12$ to end up in an $L_2$-space, simply
because $d^{a-1}$ is only in $L_2(\Omega )$ then; the factor $d^{a-1}$ is provided by
$K^{a-1}_{(0)}$ no matter how large $s$ may be. Recall also the
observation in (5.5) that for $a>\frac12$, $H^{(a-1)(2a)}(\comega
)\subset L_2(\Omega )$, and the description in Example 5.5 of how factors $d^{a-1}$ come in. Then we can show the following theorem:

\proclaim{Theorem 6.3} Assume Hypothesis {\rm 3.1}, For $a>\frac12$
and $\tau >2a+1$, consider the evolution problem
$$\aligned
Pu+\partial_tu&=f\text{ on }\Omega \times I ,\\
u&=0\text{ on }(\Bbb R^n\setminus\Omega )\times I, \\
\gamma _0^{a-1}u&=\psi \text{ on }\partial\Omega \times I,\\
u|_{t=0}&=0.
\endaligned \tag6.13
$$

For $f(x,t)$ given in $L_2(\Omega \times I)$, and $\psi (x,t)$
  given in $L_2(I; H^{a+\frac12}(\partial\Omega ))\cap\ol H^1(I;
  H^{\varepsilon }(\partial\Omega ))$ with $\psi (x,0)=0$ (some $\varepsilon >0$), there is a
  unique solution $u(x,t)$ of {\rm (6.13)} satisfying
  $$
u\in L_2(I;H^{(a-1)(2a)}(\comega ))\cap \ol H^1(I; L_2(\Omega )).\tag6.14
$$
\endproclaim

\demo{Proof} Let $v(x,t)=K^{a-1}_{(0)}\psi (x,t)$; it lies in
$L_2(I;H^{(a-1)(2a)}(\comega ))\cap \ol
H^1(I;H^{(a-1)(a-\frac12 +\varepsilon )}(\comega ))$ in view of
(6.12), contained in $ \ol H^1(I; L_2(\Omega ))$ by Lemma 5.3. It
satisfies 
$$
\gamma _0^{a-1}v=\psi ,\quad v|_{t=0}=0,\quad r^+Pv \in L_2(\Omega
\times I),\quad \partial_tv\in L_2(\Omega \times I).
$$
Then $w=u-v$ is in $L_2(I; H^{(a-1)(2a)}(\comega))$ with $\gamma
_0^{a-1}w=0$, hence in $L_2(I; H^{a(2a)}(\comega))$ by Theorem 2.3.
Moreover, $(r^+P+\partial_t)(u-v) \in L_2(\Omega \times I)$.
Thus in order for $u$ to solve (6.13), $w$ must 
solve a problem (6.7) with homogeneous boundary
condition and $f$ replaced by $f-(r^+P+\partial_t)v$. Here Theorem 6.2 $1^\circ$ assures that there is a unique
solution $w\in L_2(I;H^{a(2a)}(\comega ))\cap \ol H^1(I; L_2(\Omega
))$. Then $u=v+w$ is the unique solution of (6.13), satisfying (6.14).  \qed
\enddemo

Also the statements in Theorem 6.2 with higher derivatives in $t$ can
be extended to nonhomogenous problems; we leave this to the interested
reader. As for a lifting of the $x$-regularity, there is very little
leeway ($0\le s< a-\frac12$), so we leave out details.

\subsubhead 6.2 General $q\in \,]1,\infty [\,$ \endsubsubhead

For general $q$, the case of translation-invariant operators with real even homogeneous
symbol has been treated in \cite{G18a,G18b} for smooth domains $\Omega
$.  We shall present some straightforward consequences for our
types of nonsmooth domains, and then supply this with new results for
nonhomogeneous boundary conditions.

\proclaim{Theorem 6.4} Assume Hypothesis {\rm 3.1}, and assume
moreover that the symbol $p$ is independent of $x$, and is real, even
and homogeneous of degree $2a$.

$1^\circ$ The evolution problem {\rm (6.7)} for $P$ with homogeneous
Dirichlet condition has for every $f\in
L_q(\Omega \times I)$ a unique solution
  $$
u\in L_q(I;H_q^{a(2a)}(\comega ))\cap \ol H_q^1(I; L_q(\Omega ));\tag6.15
$$
moreover,
$$
u\in \ol C^0(I;L_q(\Omega )).\tag 6.16
$$

$2^\circ$ Let $s\in \rp\setminus \N$. Then when $u$ solves {\rm
(6.7)} with $I$ replaced by $\rp$, 
$$
f\in \dot C^s(\crp ; L_p(\Omega ))\iff u\in 
\dot C^s(\crp ;
H_q^{a(2a)}(\comega ))\cap 
\dot C^{s+1}(\crp ;
L_q(\Omega )). \tag6.17
$$

\endproclaim

\demo{Proof} $1^\circ$ was proved for $C^\infty $-domains in
\cite{G18a, Th.\ 4.3}, also recalled in \cite{G18b, Th.\ 5.9}. Since
the symbol is $x$-independent, the only new aspect in the nonsmooth
case is that $\Omega $ is allowed to be $C^{1+\tau }$. The details
of proof given in \cite{G18a} are still valid in that case: Let
$k(y)=\Cal F^{-1}p(\xi )$, it is homogeneous of degree $-2a-n$,
$C^\infty $, positive and even for $y\ne 0$, and the sesquilinear form
defining $P_{D,2}$ can be written as
$$
Q(u,v)=\tfrac12\int_{\rn }(u(x)-u(y))(\bar v(x)-\bar
v(y))k(x-y)\,dxdy\text{ for }u,v\in \dot H^a(\comega)
$$
(cf.\ also Ros-Oton \cite{R16}). As noted in Example 4.20, the operator $P_{D,2}$ it defines, as well
as the operators $P_{D,q}$, are bijective for all $1<q<\infty $. 
Moreover, the quadratic form $E(u)=Q(u,u)$ with domain $D(E)=\dot
H^a(\comega)$ has the Markovian property:
When $u_0$ is defined from a real
function $u\in D(E)$ by $
u_0=\min\{\max\{u,0\},1\}$, then $u_0\in D(E)$ and $E(u_0)\le
E(u)$. It is a so-called Dirichlet form, as explained in  Fukushima, Oshima and Takeda
\cite{FOT94}, pages 4--5 and  Example 1.2.1, and Davies
\cite{D89}. (Such a Markovian property enters e.g.\ in \cite{BBC03} for the
regional Dirichlet problem.)
 Then, by \cite{FOT94} Th.\ 1.4.1 and \cite{D89} Th.\ 1.4.1--1.4.2,
 $-P_{D,q}$
generates a strongly continuous contraction semigroup $T_q(t)$ not only in $L_2(\Omega )$ for $q=2$ but also
in $L_q(\Omega )$ for any $1<q<\infty $, and $T_q(t)$ is bounded holomorphic.
Hereby we have the prerequisites to apply the theorem of Lamberton
\cite{L87}, which shows that
(6.7) is solvable with $u\in L_q(I;D(P_{D,q}))\cap \ol H_q^1(I;
L_q(\Omega ))$.
Now we know moreover from (4.14) (based on Theorem 3.2) that
$D(P_{D,q})=H_q^{a(2a)}(\comega)$, so $1^\circ$ follows by insertion
of this fact. Since $D(P_{D,q})\subset L_q(\Omega )$, the statement
$u\in \ol C^0(I;L_q(\Omega ))$ follows from the continuity in $t\in
\ol I$ of
functions in $\ol H_q^1(I; X)$ valued in a Banach space $X$.

$2^\circ$. The details for this extension were given in 
\cite{G18b, Sect.\ 5.3}, where it is proved (by use of Hille and
Phillips \cite{HP57, Th.\ 17.5.1}, also in Kato  \cite{K66,
Th.\ IX.1.23}) that the operator properties shown
for $P_{D,q}$ moreover imply a resolvent estimate
$$
\ang\lambda \|(P_{D,q}-\lambda )^{-1}\|_{\Cal L(L_q(\Omega ))}\le C,\tag6.18
$$
for $\lambda $ outside a sectorial region like $M$ (4.9). Then (6.17)
follows as in \cite{G18b,Th.\ 5.14}
by use of a theorem of Amann \cite{A97} cited as  \cite{G18b,Th.\ 5.13}. \qed 
\enddemo

The analysis shows in particular that $P_{D,q}$ has {\it maximal
$L_q$-regularity in $I$} as defined e.g.\ in Denk and Seiler
\cite{DS15}. Also the other results of \cite{G18b, Sect.\ 5.3} extend
to nonsmooth $\Omega $, with $I=\rp$
or $\,]0,T[\,$. 

We can now moreover show results for the problem (6.13) with
a nonhomogeneous local Dirichlet condition, when $q<(1-a)^{-1}$:

\proclaim{Theorem 6.5} Assumptions of Theorem {\rm 6.4}. If in addition $\tau >2a+1$ and $q<(1-a)^{-1}$, the evolution problem {\rm (6.13)} for $P$ with nonhomogeneous
Dirichlet condition has for every $f\in
L_q(\Omega \times I)$, $\psi \in L_q(I;
B_q^{a+1/{q'}}(\partial\Omega ))\cap \ol H_q^1(I; B_q^{\varepsilon
}(\partial\Omega ))$ with $\psi |_{t=0}=0$ (some $\varepsilon >0$), a unique solution
  $$
u\in L_q(I;H_q^{(a-1)(2a)}(\comega ))\cap \ol H_q^1(I; L_q(\Omega )).\tag6.19
$$
\endproclaim

\demo{Proof} One proceeds exactly as in the proof of Theorem 6.3,
eliminating the boundary condition by subtracting a lifting of the
boundary value (using Theorem 2.3), and applying the result for the
homogeneous case. \qed

\enddemo

Consequences can be drawn as in Theorem 6.4 $2^\circ$ concerning
higher time-derivatives.

 Hereby all values of $a\in
\,]0,1[\,$ can be included in the treatment of evolution problems with
nonhomogeneous boundary
conditions.

Theorems 6.4 and 6.5 apply to $P=(-\Delta )^a$, and to all operators with
symbols $p(\xi )=g(\xi /|\xi |)|\xi |^{2a}$,  where $g(\eta )$ is
a real positive even $C^\infty $-function on $S^{n-1}$.
(Note the explanation after (2.8) about how the nonsmoothness at $\xi =0$
can be handled.)

\subhead A. Appendix \endsubhead

\subsubhead A.1 Identification of weighted boundary
maps \endsubsubhead

In \cite{A15}, Abatangelo showed an integral formula, that
we  recall here  for functions $u$
vanishing on $\rn\setminus \comega $, $u$ and $v$ being real: 

\proclaim{Proposition A.1}\cite{A15}
Let $\Omega $ be open, bounded and $C^{1,1}$.
Let $u$ be a function such that $u\in C^{2a+\varepsilon }_{loc}(\Omega)$, $u/d^{a-1}\in C(\comega)$
and $u=0$ on $\rn\setminus\comega $. Let $v\in C^a(\rn)$ be such that $v=0$ on
$\rn\setminus\comega $ and $r^+(-\Delta )^av\in C_0^\infty (\Omega ) $.

Define a boundary value of $u$ in the following way:
With $G_\Omega (x,y)$ denoting the Green's kernel for the restricted fractional
Laplacian on $\Omega $ with homogeneous Dirichlet condition (i.e., the kernel of
the operator solving {\rm (1.7)} with $P=(-\Delta )^a$), let
$$
M_\Omega (x,\theta )=\lim_{y\to\theta , y\in\Omega
}\frac{G_\Omega (x,y)}{d(y)^a}\text{ for }x\in\Omega ,\theta \in\partial\Omega
,\quad
m(x)=\int_{\partial\Omega }M_\Omega (x,\theta ')d\sigma (\theta ') ,
$$
and define the trace operator $E$ by
$$
Eu(\theta )=\lim_{x\to\theta ,x\in\Omega
}\frac{u(x)}{d^{a-1}(x)m(x)},\quad \theta \in\cdot\Omega .
$$
Then (cf.\ \cite{A15, (9)}): 
$$
\int_\Omega u\,(-\Delta )^{a}v\,dx-\int_\Omega (-\Delta
)^{a}u\,v\,dx=\int_{\partial\Omega }Eu\,\gamma _0(\tfrac v{d^a})\, d\sigma .
\tag A.1
$$
\endproclaim

It is stated in \cite{A15, Sect.\ 3.1}, that $Eu\in
C(\partial\Omega )$.

In \cite{G18} we showed a Green's formula, which implies the
following ``halfways  Green's formula'' (cf.\ Cor.\ 4.5 (4.34) there):

\proclaim{Proposition A.2}\cite{G18}
Let $\Omega \subset\rn$, open, bounded  and smooth, and let $P$  be a
classical 
pseudodifferential operator of order $2a>0$  satisfying
the $a$-transmission condition at $\partial\Omega $ (it suffices for
this that $P$ is {\it even}). If $u\in H^{(a-1)(s)}(\comega)$, $v\in
H^{a(s)}(\comega)$,  $s>a+\frac12$, and $s\ge 2a$, 
$$
\int_\Omega u\,\overline{P^*v}\,dx-\int_\Omega Pu\,\bar v\,dx
=\Gamma (a)\Gamma (a+1)\int_{\partial\Omega }s_0\gamma _0(\tfrac{u}{d^{a-1}})\gamma
_0(\tfrac{\bar v}{d^{a}})
\, d\sigma .
\tag A.2
$$
\endproclaim

Here $s_0(x)$ is a function defined from the principal symbol of $P$; it is
1 when $P=(-\Delta )^a$. The spaces $ H^{(a-1)(s)}(\comega)$ and
$H^{a(s)}(\comega)$ are the solution spaces for the nonhomogeneous, 
resp.\ homogeneous, Dirichlet problem for $P$ in the scale of
$L_2$-Sobolev spaces.

The formula holds a fortiori for functions in H\"older spaces
$$
\aligned
u&\in C_*^{(a-1)(s+\varepsilon )}(\comega)\subset \dot
C^{s+\varepsilon }(\comega)+d^{a-1}e^+C^{s-a+1+\varepsilon
}(\comega),\\
v&\in C_*^{a(s+\varepsilon )}(\comega)\subset \dot C^{s+\varepsilon
}(\comega)+d^{a}e^+C^{s-a+\varepsilon }(\comega)\subset
d^{a}e^+C^{s-a+\varepsilon }(\comega);
\endaligned \tag A.3
$$
here $d^{a-1}$ is $L_1$-integrable over $\Omega $, and we take small
$\varepsilon >0$ such that the indexations avoid integers.

The formula (A.2) can be extended to suitable nonsmooth domains by use of
the solvability results in Section 3 above, but we shall not take up space
here with details.

Applying (A.1) and (A.2) to  $P=(-\Delta )^a=P^*$, we find that
$$
\int_{\partial\Omega }Eu\,\gamma _0(\tfrac v{d^a})\, d\sigma =\Gamma (a)\Gamma (a+1)\int_{\partial\Omega }\gamma _0(\tfrac{u}{d^{a-1}})\gamma
_0(\tfrac{ v}{d^{a}})
\, d\sigma \tag A.4
$$
holds for real functions $u$ and $v$ satisfying the hypotheses for both
propositions.

In view of (A.3) and the fact that $C_*^{(a-1)(2a+\varepsilon )}
(\comega)\subset C^{2a+\varepsilon }_{\operatorname{loc}}(\Omega )$,
the
requirements on $u$ in Proposition A.1 are satisfied when
$u\in C_*^{(a-1)(2a+\varepsilon )}(\comega)$.
As for $v$,  we know from \cite{G15a} that $v\in \dot C^a(\comega)$ with $(-\Delta )^av\in C^\infty
(\comega)$ implies that  $v\in \Cal E_a(\comega)$ (cf.\ (1.1)), and $v\in \Cal
E_a(\comega)$ implies $r^+(-\Delta )^av\in C^\infty (\comega)$, so the functions $v$ satisfying the requirement of
Proposition A.1 are a subset of $\E_a(\comega)$,
We claim that one can conclude
$$
Eu=c\gamma _0(\tfrac u{d^{a-1}}),\quad c=\Gamma (a)\Gamma (a+1).
\tag A.5
$$
To show this, we still have to prove that $\gamma _0(v/d^a)$ assumes
enough values to allow a passage from the weak identity (A.4) to the
identity (A.5).

\proclaim{Lemma A.3} For the functions $v$ satifying the hypotheses of Proposition {\rm
A.1}, $\gamma _0(v/d^a)$ runs through a dense subset of $L_2(\partial\Omega )$.
\endproclaim

\demo{Proof}
For small $\varepsilon '>0$, $\ol H^{\frac12-\varepsilon '}(\Omega )$ identifies with
$\dot H^{\frac12-\varepsilon '}(\comega)$,
and $C_0^\infty (\Omega )$ is dense in the space. The solution
operator for the homogeneous Dirichlet problem maps this space
bijectively onto
$H^{a(2a+\frac12-\varepsilon ')}(\comega)$. So when $(-\Delta )^av$ runs through
the dense subset $C_0^\infty (\Omega )$ of  $\ol
H^{\frac12-\varepsilon '}(\Omega )$ ,  $v$ runs through a dense
subset of $ H^{a(2a+\frac12-\varepsilon ')}(\comega)$.

Here we have from \cite{G19}
Th.\ 3.4 the precise statement (defining a direct sum) when
$2a+\frac12-\varepsilon '-a\in \,]\frac12,\frac32[\,$, i.e., $0<a-\varepsilon '<1$:
$$
v\in H^{a(2a+\frac12-\varepsilon ')}(\comega)\iff
v=w+d^aK_{(0)}\varphi ,\;w\in \dot
H^{2a+\frac12-\varepsilon '}(\comega), \;
\varphi \in H^{a-\varepsilon '}(\partial\Omega),
$$
where  $\varphi =\gamma _0(v/d^a)$, and $K_{(0)}$ is a certain Poisson
operator with $\gamma _0K_{(0)}=I$ (this was first shown in local coordinates
in \cite{G15a}, Th.\ 5.4, that we could point to instead). Thus when $v$ runs through a dense
subset of $ H^{a(2a+\frac12-\varepsilon ')}(\comega)$, $\varphi $ runs
through a dense subset of $H^{a-\varepsilon '}(\partial\Omega)$.
A fortiori, for the considered
$v$,  $\varphi =\gamma _0(v/d^a)$ runs through a dense subset of
$L_2(\partial\Omega ) $.\qed

\enddemo

Then the conclusion of (A.5) from (A.4) follows, and we have obtained:

\proclaim{Theorem A.4} When $\Omega $ is bounded smooth,
the boundary trace $Eu$ introduced in \cite{A15} equals the constant $c=\Gamma (a)\Gamma (a+1)$ times
$\gamma _0(u/d^{a-1})$, cf.\ {\rm (A.5)}.

This holds for functions $u\in
C^{2a+\varepsilon }_{\operatorname{loc}}(\Omega )$ satisfying {\rm (A.3)}. More
precisely,  they are the solutions of the nonhomogeneous Dirichlet problem
{\rm (1.10)}
with $f\in C^\varepsilon (\comega)$,
$\varphi \in C^{a+\varepsilon }(\partial\Omega ) $.

\endproclaim    

A different proof is given in \cite{CGV21}, referring to the Pohozaev
formula shown by Ros-Oton and Serra in \cite{RS14}. In \cite{CGV21},
the constant 
$c$ is stated to be equal to \linebreak $\Gamma (a+1)^2$ (same constant as in
the Pohozaev formula); it seems that a factor $a^{-1}$ has been overlooked
in the application of the boundary mapping to a derivative of
$u$. (This is confirmed in an arXiv posting 22.4.2022, arXiv:2004.04579v2.)

\subsubhead A.2 An embedding property in H\"older
spaces \endsubsubhead

For a $C^{1+\tau }$-domain, the distance function $d_0(x)$, equal to
$\operatorname{dist}(x,\partial\Omega )$ near $\partial\Omega $ and extended
smoothly and positively to $\Omega $, is a $C^\tau $-function on
$\comega $. If  $\tau \ge 1$, $d_0$ is $C^{\tau +1}$ (as recalled in
Section 2).

\proclaim{Lemma A.5} Let $\Omega $ be a $C^{1+\tau }$-domain, $\tau
>0$.
Let $a$ and $b>0$ with $a+b<1+\tau $. For $a,b,a+b\notin{\Bbb N}$,
there holds
$$
\dot C^{a+b}(\comega)\subset d_0(x)^a \dot C^{b
}(\comega ).\tag A.6
$$

\endproclaim

\demo{Proof} The following calculations take place in a neighborhood
of $\partial\Omega $ where $d_0(x)=\operatorname{dist}(x,\partial\Omega )
$. On the interior, $v=u/d_0^a$ is $C^{b}$ simply because $u$ and
$d_0^{-a}$ are so.

$1^\circ$. First consider the basic case where $a+b<1$.
Let $u\in \dot C^{a+b}(\comega)$, and let
$v(x)=u(x)/d_0(x)^a$. 
For $u$ we have
$$
|u(x)-u(y)|\le C|x-y|^{a+b},\quad|u(x)|\le C d_0(x)^{a+b}. \tag A.7
$$

For the various positions of $x$ and $y$, we distinguish the
following two cases: a) $|x-y|>\frac13 d_0(y)$, b) $|x-y|\le
\frac13 d_0(y)$. We can assume $0<d_0(x)\le d_0(y)$ (and of course $x\ne y$).
Denote by $x_0$ a point on $\partial\Omega
$ where $|x-x_0|=d_0(x)$; then $$
\aligned
&d_0(y)\le  
|y-x_0|\le
|y-x|+|x-x_0|=|y-x|+d_0(x),\\
&\text{ hence }|d_0(x)-d_0(y)|\le |x-y|.
\endaligned\tag A.8
$$

a) $|x-y|>\frac13 d_0(y)$. Here
$$
\aligned
\frac{|v(x)-v(y)|}{|x-y|^b}&=\frac{|u(x)/d_0(x)^a-u(y)/d_0(y)^a|}{|x-y|^b}
\le
\frac{Cd_0(x)^{b}}{|x-y|^{b}}+\frac{Cd_0(y)^{b}}{|x-y|^{b}}\\
&\le C d_0(x)^b(\tfrac13 d_0(x))^{-b}+C d_0(y)^b(\tfrac13 d_0(y))^{-b}=C'.
\endaligned\tag A.9
$$

b) $|x-y|\le\frac13 d_0(y)$. In view of (A.8)  we  have 
$$  
d_0(x)\ge d_0(y)-|x-y|\ge \tfrac23 d_0(y)\ge 2|x-y|.
$$
Now
$$
\aligned
\frac{|v(x)-v(y)|}{|x-y|^b}&=\frac{|u(x)/d_0(x)^a-u(y)/d_0(x)^a+u(y)/d_0(x)^a-u(y)/d_0(y)^a|}{|x-y|^b}\le
I+II,\\
I&=\frac{|u(x)-u(y)|d_0(x)^{-a}}{|x-y|^b},\quad
II=\frac{|u(y)||d_0(x)^{-a}-d_0(y)^{-a}|}{|x-y|^b}.
\endaligned\tag A.10
$$
For $I$,
$$
I\le C|x-y|^ad_0(x)^{-a}\le C(\tfrac12 d_0(x))^a(d_0(x))^{-a}= C''.
$$
For $II$, denote $d_0(x)=s$, $d_0(y)=t$ and $r=1-s/t$. If $d_0(x)=d_0(y)$, 
$II=0$, so we can assume $s\ne t$. Here
$$
0< r=\frac {t-s}{t}\le \frac{|x-y|}{3|x-y|}=\tfrac14.
$$
Then in view of (A.7) and (A.8),
$$
\aligned
II&\le
Ct^{a+b}\frac{|s^{-a}-t^{-a}|}{|s-t|^b}=C\frac{|1-(s/t)^{-a}|}{|1-s/t|^b}=
C\frac{|1-(1-r)^{-a}|}{r^b}\\
&=C\frac{|ar+O(r^2)|}{r^b}\le C''',
\endaligned \tag A.11
$$
since $r\in \,]0,\frac13]$ and $b<1$ (we have used a  Taylor expansion of $(1-r)^{-a}$).
This shows that $v\in C^b(\comega)$. Since $v(x)$ is $O(d_0(x)^b)$ with $b>0$, $v$ vanishes at
$\partial\Omega $, so in fact, $v\in \dot C^b(\comega)$.

$2^\circ$. Next, consider cases where $a$ is let free, but $b$ is
still assumed to be $<1$. So $a+b\in \,]k,k+1[\,$, where the case $k=0$
was treated above. Let $k=1$. In (A.7), the estimates are replaced by
$$
\aligned
|u(x)-u(y)|&\le C|x-y|,\;|\nabla u(x)-\nabla u(y)|\le
 C|x-y|^{a+b-1},\\
  |u(x)|&\le Cd_0(x)^{a+b},\; |\nabla u(x)|\le Cd_0(x)^{a+b-1}.
  \endaligned
  \tag A.12
$$
The calculation (A.9) carries over
verbatim, and so does the treatment of $II$ defined in (A.10). Only $I$
needs a modified argument: There holds
for $\theta \in [0,1]$:
$$
|d_0(x+\theta (y-x))-d_0(x)|\le |x+\theta (y-x)-x|\le |x-y|,
$$
hence  since $|x-y|\le \frac12 d_0(x)$, 
$$
\aligned
|u(x)-u(y)|&=|(x-y)\cdot\int_0^1\nabla u(x+\theta (y-x))\,d\theta |\le
C|x-y|\sup_\theta (d_0(x+\theta (y-x))^{a+b-1}\\
&\le
C|x-y|(d_0(x)+|x-y|)^{a+b-1}\le C'|x-y|d_0(x)^{a+b-1}.
\endaligned
$$
Thus 
$$
I\le C'|x-y|d_0(x)^{a+b-1}d_0(x)^{-a}|x-y|^{-b}=C'|x-y|^{1-b}d_0(x)^{b-1}\le C'(\tfrac12)^{b-1}.
$$
For higher $k$, one similarly uses Taylor's formula on the
$k$'th level.

$3^\circ$. Finally, let also $b>1$. In this case, $\tau >1$, so
$d_0(x)\in C^{1+\tau }(\comega)$. When $b\in \,]k,k+1[\,$, it is the
$k$'th derivatives $\partial^\alpha v$ with $|\alpha |=k$
that we have to estimate. E.g., for $k=1$,
$$
\partial_jv=\partial_j(ud^{-a})=\partial_ju\,d^{-a}+u(-a)d^{-a-1}\partial_jd.
$$
For $\partial_ju\,d^{-a}$, the result follows by application of $2^\circ$ to $\partial_ju\in
\dot C^{a+b-1}$, this gives an element of $\dot C^{b-1}$. For $u\,d^{-a-1}\partial_jd$, we
likewise apply $2^\circ$ to find an element of $\dot C^{b-1}(\comega)$, since higher
positive $a$ are allowed, and multiplication by $\partial_jd\in
C^{\tau }(\comega)$ preserves being in $\dot C^{b-1}(\comega)$. Altogether,
$\partial_jv\in \dot C^{b-1}(\comega)$, and hence $v\in
\dot C^b(\comega)$.

Passing to $k=2$, we use the information from the cases $k=0$ and
1. Then we go on successively to all the relevant larger $k$; by the Leibniz
formula there are more and more elements to account for, but the
principle is the same.\qed

\enddemo

The lemma can be extended to the distance function $d(x)=x_n-\zeta (x')$
defined for the set $\R^n_\zeta =\{x=(x',x_n)\mid x_n>\zeta (x')\}$,
$\zeta \in C^{1+\tau }(\R^{n-1})$;
this is useful when $\tau <1$, since the regularity results in the
main text are
primarily formulated with  $d$ rather than $d_0$ then.

\proclaim{Lemma A.6} Let  $\tau >0$ and consider the  curved halfspace  $\R^n_\zeta
=\{x=(x',x_n)\mid x_n>\zeta (x')\}$, defined from a function
$\zeta \in C^{1+\tau }(\R^{n-1})$, and provided with the distance
function $d(x)=x_n-\zeta (x')$ for $x_n\le K+1$ (some  $K\ge \sup |\zeta |$).

Let  $a$ and $b>0$ with $a+b<1+\tau $. For $a,b,a+b\notin{\Bbb N}$, one has near $\partial\R^n_\xi $ that
$$
u\in \dot C^{a+b}(\R^n_\zeta )\implies u(x)/ d(x)^a\in \dot C^{b
}(\R^n_\zeta  ).\tag A.13
$$
When $b<\tau $, the same result holds if $d $ is replaced by another distance
function $d '$ in $C^{1+\tau }(\rn_\zeta )$ bounded above and below by $d$.
\endproclaim

\demo{Proof} We shall show how the proof of Lemma A.5 is generalized
to this case. Assume $a+b<1$. Let $u\in \dot C^{a+b}(\R^n_\zeta )$, and let
$v(x)=u(x)/d(x)^a$ (where $d(x',x_n)=x_n-\zeta (x')>0$). 
For $u$ we have
$$
|u(x)-u(y)|\le C|x-y|^{a+b},\quad|u(x)|\le C d(x)^{a+b}. \tag A.14
$$
Let $x$ and $y$ be points with $d(y)>d(x)>0$. Since $\zeta \in  C^{1+\tau
}(\R^{n-1})$, there is a $c_1$ such that $|\zeta (y')-\zeta
(x')|\le c_1 |y'-x'|$. For the given
$x=(x',x_n)$, denote $(x',\zeta (x'))=x_0$, with a similar notation for $y$.
$x_0$ is the  ``footpoint''   on $\partial\R^n_\zeta $, i.e.\ the point with
the same $x'$-value as $x$, so that $d(x)=|x-x_0|$. 

We then have
$$
|y-x_0|\le |y-x|+|x-x_0|=|y-x|+d(x).
$$
Moreover, 
$$
|y-x_0|=|y-y_0+y_0-x_0|\ge d(y)-|y_0-x_0|\ge d(y)-(1+c_1 )|x'-y'|. 
$$
Then if we assume $|x-y|\le \frac13 (2+c_1)^{-1}d(y)$, then 
$$
d(x)\ge |y-x_0|-|y-x|\ge d(y)-(2+c_1)|y-x|\ge \tfrac23 d(y).\tag A.15
$$


Now proceed as in the proof of Lemma A.5:

a) The case $|x-y|>\frac13 (2+c_1)^{-1}d(y)$. Here the desired estimate is
obtained as in a) of Lemma A.5.

b) The case $|x-y|\le \frac13 (2+c_1)^{-1}d(y)$. In view of (A.15) we have 
$$
d(x)\ge \tfrac23 d(y)\ge \tfrac23 3(2+c_1)|x-y|.
$$
With these inequalities, the proof of b) in Lemma A.5 goes through, just with other
constants.

In this way the desired result is obtained for $a+b<1$; and it extends
to more general values in the same way as in Lemma A.5.

The last statement follows from the fact that $d/d'$ and $d'/d$ are
$C^\tau $; it is included to allow general
distance functions as introduced around (2.2).
\qed

\enddemo

\subhead Acknowledgement \endsubhead

The author is grateful to Helmut Abels and Xavier Ros-Oton for useful conversations.

\enddocument

\enddocument